\documentclass[10pt]{article}

\usepackage{amssymb,amsmath,amsthm}
\usepackage{caption,graphicx}
\usepackage{paralist}
\usepackage[text={6.75in,9.5in},centering,letterpaper]{geometry}
\usepackage[comma,square,sort&compress,numbers]{natbib}

\setcaptionmargin{0.25in}

\newtheorem{theorem}{Theorem}[section]
\newtheorem{lemma}[theorem]{Lemma}
\newtheorem{corollary}[theorem]{Corollary}
\newtheorem{proposition}[theorem]{Proposition}
\newtheorem{hypothesis}{Hypothesis}

\makeatletter\@addtoreset{equation}{section}\makeatother

\begin{document}

\title{Localized radial roll patterns in higher space dimensions}

\author{
Jason J. Bramburger\thanks{Division of Applied Mathematics, Brown University} \and
Dylan Altschuler\thanks{Department of Mathematics, New York University} \and
Chloe I.~Avery\thanks{Department of Mathematics, University of Chicago} \and
Tharathep Sangsawang\thanks{Department of Mathematics, University of Texas at Austin} \and
Margaret Beck\thanks{Department of Mathematics, Boston University} \and
Paul Carter\thanks{Department of Mathematics, University of Arizona} \and
Bj\"orn Sandstede\footnotemark[1]
}

\date{}
\maketitle

\begin{abstract}
Localized roll patterns are structures that exhibit a spatially periodic profile in their center. When following such patterns in a system parameter in one space dimension, the length of the spatial interval over which these patterns resemble a periodic profile stays either bounded, in which case branches form closed bounded curves (``isolas"), or the length increases to infinity so that branches are unbounded in function space (``snaking"). In two space dimensions, numerical computations show that branches of localized rolls exhibit a more complicated structure in which both isolas and snaking occur. In this paper, we analyse the structure of branches of localized radial roll solutions in dimension 1+$\varepsilon$, with $0<\varepsilon\ll1$, through a perturbation analysis. Our analysis sheds light on some of the features visible in the planar case.
\end{abstract}


\section{Introduction}

\begin{figure}[t]
\centering
\includegraphics{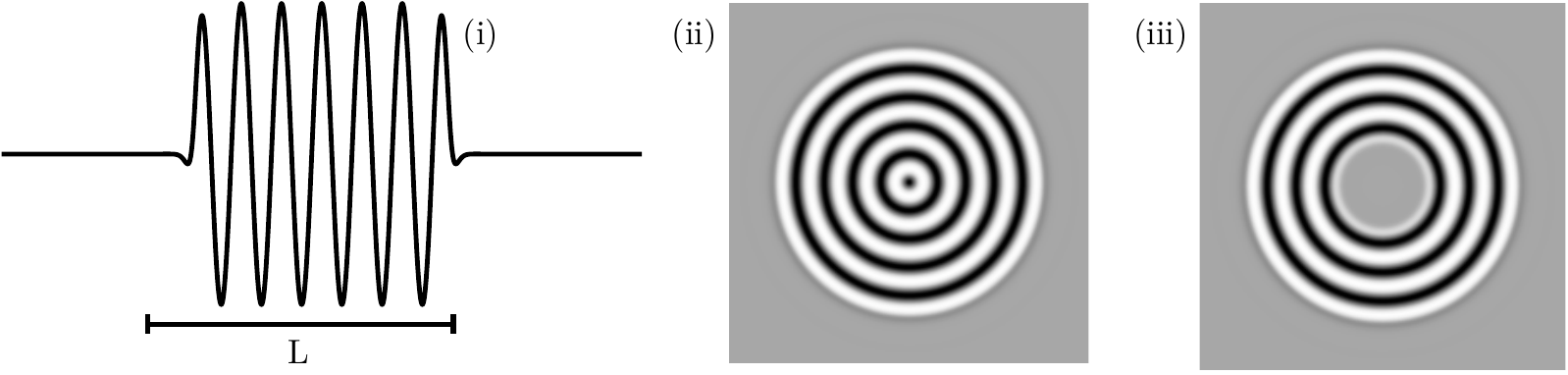}
\caption{Panel (i) shows the profile of a localized roll pattern with plateau length $L$ as a function of the spatial variable $x\in\mathbb{R}$. Panels~(ii) and~(iii) contain contour plots of planar radial spots and rings, respectively.}
\label{f:1}
\end{figure}

Spatially localized patterns can be observed in the natural world in a variety of places, such as vegetation patterns \cite{Veg1,Veg2}, crime hotspots \cite{Hotspot}, and ferrofluids \cite{Ferrofluid}. We are particularly interested in localized roll solutions. When the spatial variable $x$ is in $\mathbb{R}$, these structures are spatially periodic for $x$ in a bounded region, and they decay exponentially fast to zero as $x\to\pm\infty$; see Figure~\ref{f:1}(i) for an illustration. In planar systems with $x\in\mathbb{R}^2$, localized roll solutions may take the form of radial patterns, which are often referred to as spots and rings depending on whether the roll structures extend into the center of the pattern (spots) or not (rings); see Figure~\ref{f:1}(ii)-(iii). We refer to the length or radius of the region occupied by the periodic rolls as the plateau length of the underlying localized roll pattern.

\begin{figure}
\centering
\includegraphics{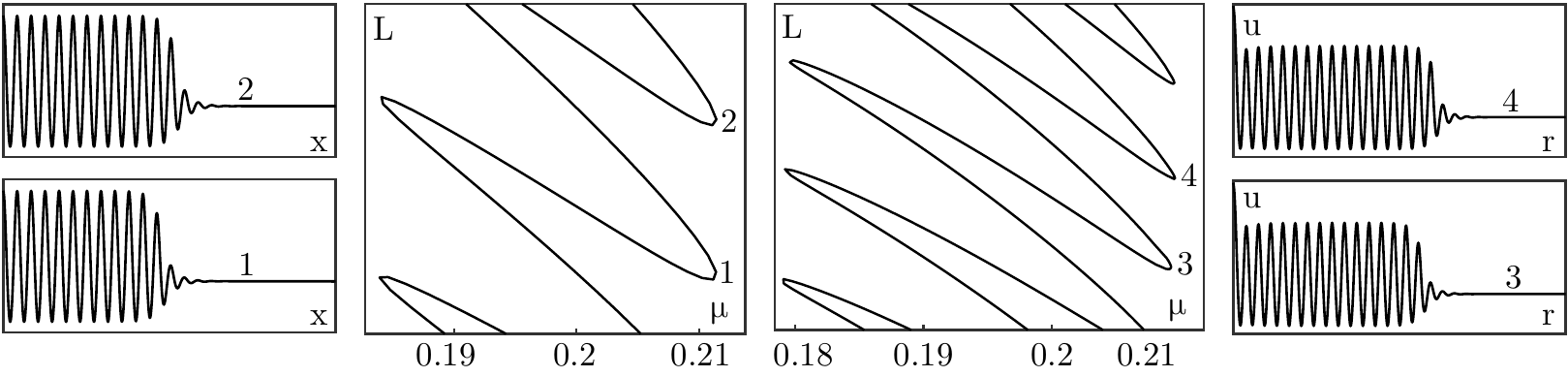}
\caption{The two center panels show the bifurcation diagrams of localized roll patterns for $n=1$ (left) and $n=2$ (right) in $(\mu,L)$-space, where $L$ denotes the plateau length as measured by the squared $L^2$-norm. The left- and rightmost panels show representative solution profiles as functions of $x$ for $n=1$ (left) and the radius $r$ for $n=2$ (right). As $L$ increases, more rolls are added to each pattern. The computations are done for (\ref{SwiftEqn}) with $\nu=1.6$.}
\label{f:2}
\end{figure}

We are interested in understanding how localized roll patterns and their plateau lengths depend on parameters. To outline the specific questions we wish to address, we focus initially on the Swift--Hohenberg equation
\begin{equation} \label{SwiftEqn}
U_t = -(1+\Delta)^2 U - \mu U + \nu U^2 - U^3, \qquad x\in\mathbb{R}^n, \quad U\in\mathbb{R},
\end{equation}
where $\Delta$ denotes the Laplace operator, $\nu$ will be held fixed, and $\mu$ is a parameter that we will vary. The Swift--Hohenberg equation admits stationary localized roll profiles in one and two space dimensions \cite{Pomeau,Coullet,Woods,Burke,Lloyd,Lloyd2,McCalla,McCalla2}. In particular, it was shown in \cite{Lloyd,McCalla2} that (\ref{SwiftEqn}) has two different stationary spot and ring patterns (referred to as spot~A, spot~B, ring~A, and ring~B) near $\mu=0$. Figure~\ref{f:2} visualizes solution branches associated with localized roll patterns of the Swift--Hohenberg equation by plotting the parameter $\mu$ for which a roll pattern exists against its plateau length $L$. As shown there, the bifurcation branches oscillate back and forth between fold bifurcations, and the plateau length increases as additional rolls are added to the pattern as each branch is traversed.

\begin{figure}
\centering
\includegraphics{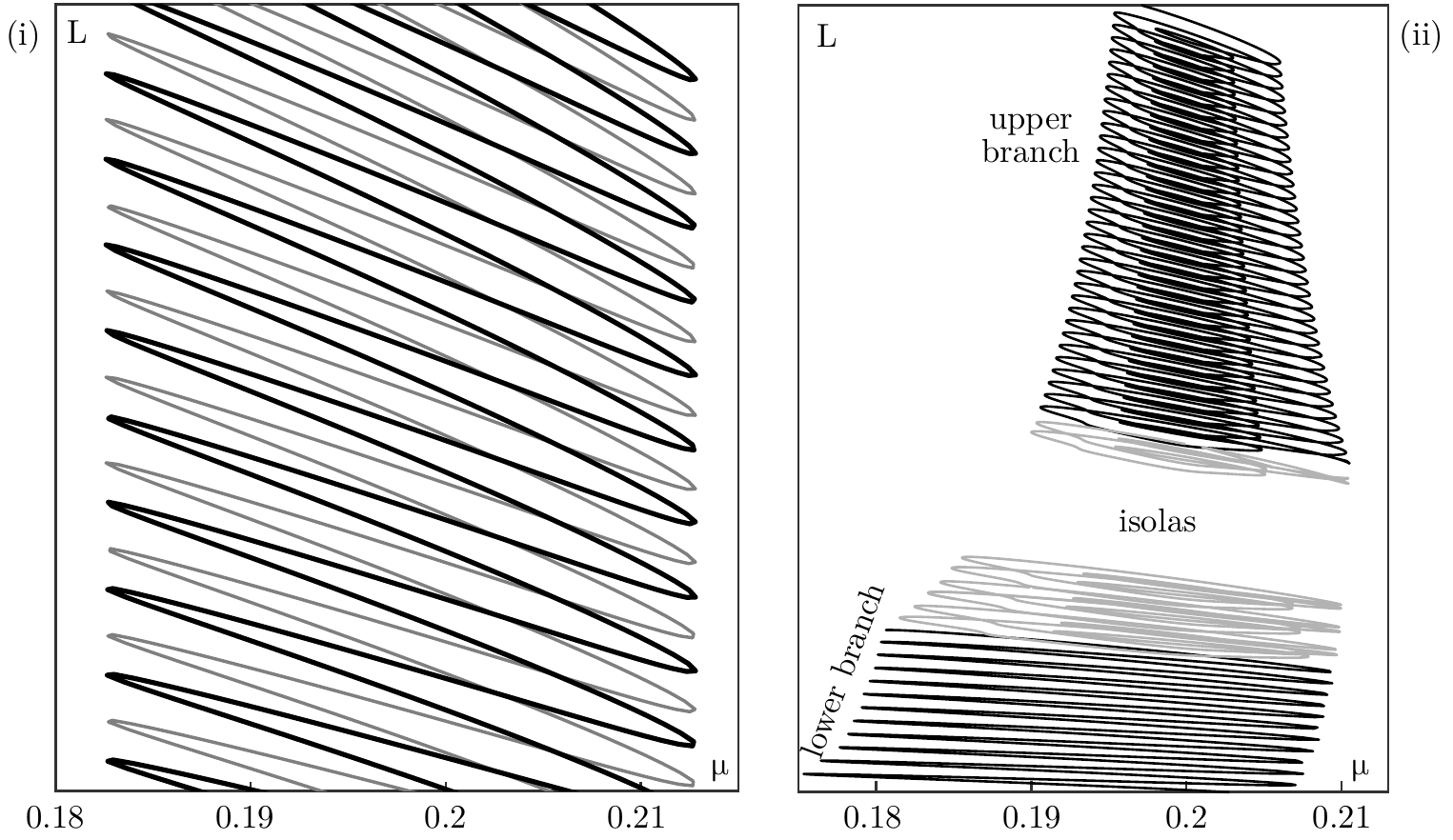}
\caption{Panels (i) and (ii) show the bifurcation diagrams of stationary localized roll patterns in dimension $n=1$ and $n=2$, respectively, of (\ref{SwiftEqn}) with $\nu=1.6$. For $n=1$, there are two solution branches, corresponding to localized rolls with, respectively, a maximum and minimum at the center, that oscillate back and forth forever. Panel~(ii) shows the solution branches corresponding to planar spot~A and ring~A patterns, which have a maximum at the center: the connected lower branch is followed by a stack of closed loops and a connected upper branch whose width shrinks as the plateau length $L$ increases. Not shown is a second similar set of branches for spot~B and ring~B solutions; we refer to \cite[Figures~2,~4, and~6]{McCalla} for the full bifurcation diagram.}
\label{f:3}
\end{figure}

A key difference between the one- and two-dimensional cases becomes apparent when the bifurcation branches are displayed over a larger range of plateau lengths. Figure~\ref{f:3}(i) shows the bifurcation diagram of one-dimensional stationary localized roll patterns: two branches exist that oscillate back and forth between two vertical aymptotes, and the profiles on these two branches differ by whether they have a minimum or a maximum at their center---we refer to these branches as snaking branches. In contrast, in the planar case, Figure~\ref{f:3}(ii) shows that the continuation of the spot and ring patterns found near $\mu=0$ leads to branches that fragment into connected lower and upper branches, which are separated by finitely many stacked closed loops that we refer to as isolas. In addition, the fold bifurcations along these branches do not align, and the width of the upper branches decreases as the plateau length increases \cite{Lloyd,McCalla,McCalla2}.

In this paper, we investigate the differences between the bifurcation diagrams in one and two space dimensions. In particular, we will analyse whether the snaking branches observed in one space dimension persist for all plateau lengths or whether they terminate at some maximal length, and we will also study whether the branch width collapses, and if so, at which value of the parameter $\mu$.

Before we outline our results, we focus briefly on the one-dimensional Swift--Hohenberg equation
\begin{equation}\label{n1}
U_t = -(1+\partial_x^2)^2 U - \mu U + \nu U^2 - U^3, \qquad x\in\mathbb{R}.
\end{equation}
For each fixed $(\mu,\nu)$, this equation admits a one-parameter family of periodic roll patterns that is parametrized by their period $p$. Intuitively, the stationary profile shown in Figure~\ref{f:1}(i) can be obtained by gluing one of these periodic profiles and the homogeneous rest state $U=0$ together. Since the steady-state equation associated with (\ref{n1}) admits the conserved quantity
\begin{equation}\label{e:H}
\mathcal{H}(U,\mu) := U_x U_{xxx} - \frac{U^2_{xx}}{2} + U^2_x + \frac{U^2}{2} + \frac{\mu U^2}{2} - \frac{\nu U^3}{3} + \frac{U^4}{4},
\end{equation}
which is conserved pointwise along each stationary solution $U(x)$ of (\ref{n1}), this quantity must vanish when evaluated along the roll pattern as $\mathcal{H}(0,\mu)=0$. This condition leads to a selection principle for the periodic profile inside a localized roll structure as there will generally be only one roll pattern $U_\mathrm{per}(x)$ in the one-parameter family for which $\mathcal{H}(U_\mathrm{per}(x),\mu)=0$. On the other hand, the Swift--Hohenberg equation is a gradient system with energy given by
\[
\mathcal{E}(U,\mu) := \int \left( \frac{(U+U_{xx})^2}{2} + \frac{\mu U^2}{2} - \frac{\nu U^3}{3} + \frac{U^4}{4} \right) \,\mathrm{d}x,
\]
and we may therefore expect that solutions with lower energy invade those with higher energy. Thus, depending on whether the energy of the selected roll pattern $U_\mathrm{per}(x)$ over one spatial period is larger or smaller than zero (the energy associated with $U=0$ vanishes), the plateau width of localized rolls should either decrease or increase as time increases. This heuristic argument shows that we may expect to observe localized roll profiles only for the single parameter value $\mu$ at which the energy $\mathcal{E}(U_\mathrm{per},\mu)$ of the selected periodic profile $U_\mathrm{per}(x)$ vanishes. This parameter value is commonly referred to as the Maxwell point $\mu_\mathrm{Max}$, and its value for $\nu=1.6$ is $\mu_\mathrm{Max}=0.2004$, which lies inside the $n=1$ snaking region shown in Figure~\ref{f:3}(i). The heuristic reason for why localized rolls exist in an open interval in parameter space, and not just at a single parameter value, is that the argument given above does not account for energy stored in the interface between the roll pattern and the homogeneous rest state. Inspecting Figure~\ref{f:3}(ii), it is tempting to conjecture that the branch in the planar case collapses onto the Maxwell point, and we will return to this conjecture below.

\begin{figure}
\centering
\includegraphics{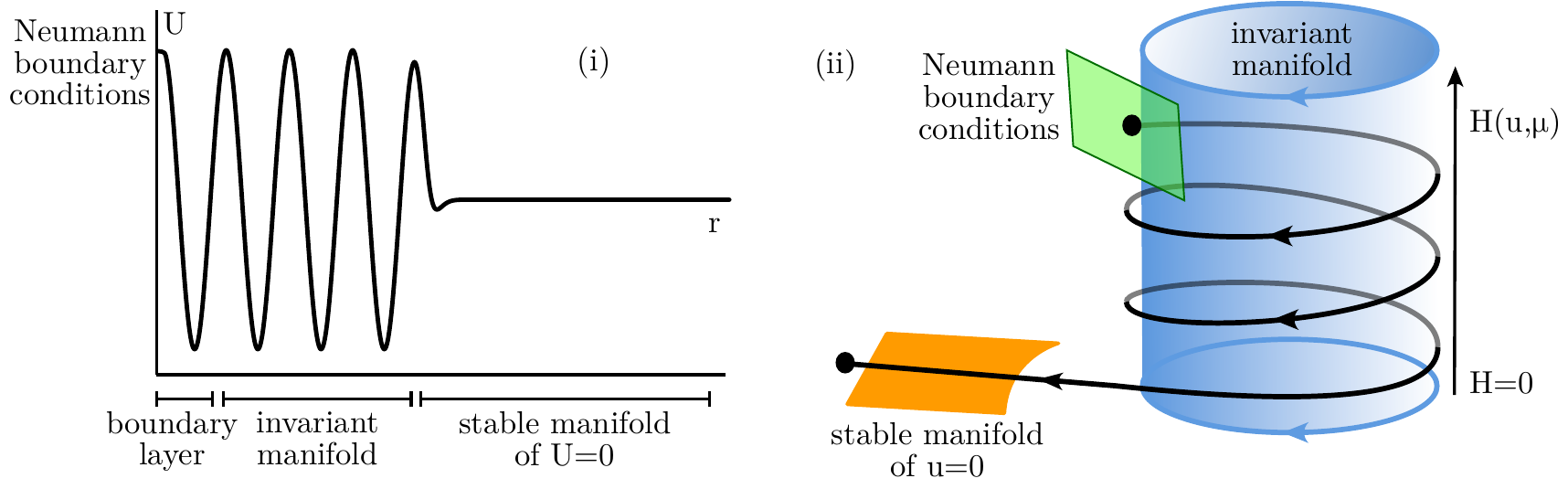}
\caption{Panels (i) and (ii) illustrate the geometry behind localized roll profiles for the radial steady-state equation (\ref{SHE_Radial}) in (i) and the first-order system (\ref{SHEqn2}) in (ii). Localized roll solutions satisfy Neumann conditions at $r=0$ and follow the invariant manifold of periodic profiles parametrized by the conserved quantity $\mathcal{H}$ before entering the stable manifold of the homogeneous rest state.}
\label{fig:Cylinder}
\end{figure}

Stationary radial solutions of the Swift--Hohenberg equation posed on $\mathbb{R}^n$ can be sought in the form $U(|x|)=U(r)$ where the profile $U(r)$ satisfies the fourth-order ordinary differential equation
\begin{equation} \label{SHE_Radial}
0 = -\left(1 + \frac{n-1}{r}\partial_r + \partial_r^2\right)^2 U - \mu U + \nu U^2 - U^3, \qquad r>0.
\end{equation}
Using the variables
\begin{equation}\label{e:nv}
u_1 = U, \quad
u_2 = U_r, \quad
u_3 = \left(1+\frac{n-1}{r}\partial_r+\partial_r^2\right) U, \quad
u_4 = \partial_r \left(1+\frac{n-1}{r}\partial_r+\partial_r^2\right) U
\end{equation}
and setting ${ }^\prime=\frac{\mathrm{d}}{\mathrm{d}r}$, we can write (\ref{SHE_Radial}) as the nonautonomous first-order system
\begin{equation}\label{SHEqn2}
\begin{split}
u_1^\prime &= u_2 \\
u_2^\prime &= u_3 - u_1 - \frac{n-1}{r}u_2 \\
u_3^\prime &= u_4 \\
u_4^\prime &= -u_3 - \mu u_1 + \nu u_1^2 - u_1^3 - \frac{n-1}{r} u_4.
\end{split}
\end{equation}
When $n=1$, equation (\ref{SHEqn2}) is autonomous and reversible under $r\mapsto-r$, and $\mathcal{H}(U,\mu)$ defined in (\ref{e:H}) continues to be a conserved quantity for (\ref{SHEqn2}) once it is rewritten in the new variables (\ref{e:nv}). The stationary periodic roll profiles of (\ref{n1}) then correspond to periodic orbits of (\ref{SHEqn2}), which form a normally hyperbolic invariant manifold $\mathcal{P}$ that is parametrized by the value of $\mathcal{H}$; see Figure~\ref{fig:Cylinder}(ii) for an illustration. To construct localized roll patterns, the approach taken in \cite{Beck} was to assume the existence of a heteroclinic orbit of (\ref{SHEqn2}) inside the invariant zero level set $\mathcal{H}^{-1}(0)$ that connects the periodic orbit in $\mathcal{H}^{-1}(0)$ to the rest state $u=0$. The analysis in \cite{Beck} then focused on constructing solutions that satisfy the Neumann boundary conditions $u_2=u_4=0$ at $r=0$ and follow the periodic orbit for $0\leq r\leq L$ with $L\gg1$ before converging to $u=0$ as $r\to\infty$: as shown in \cite{Beck}, the resulting orbits can be parametrized by their plateau length $L$ (see again Figure~\ref{fig:Cylinder}).

For $n>1$, the quantity $\mathcal{H}$ is no longer conserved for the nonautonomous system (\ref{SHEqn2}). If the perturbation terms $\mathcal{O}(|n-1|/r)$ are small, then we expect that the normally hyperbolic invariant manifold $\mathcal{P}$ persists as an integral manifold for (\ref{SHEqn2}). However, the flow on the integral manifold will no longer be periodic, and solutions may leave the cylindrical integral manifold after a finite time through its top or bottom. Thus, key to understanding the existence of localized roll patterns in higher space dimensions is to understand the dynamics on the integral manifold and to extend the analysis carried out in \cite{Beck} for the autonomous equation on $\mathcal{H}^{-1}(0)$ to the nonautonomous equation on $\mathbb{R}^4$.

Our analysis will be perturbative in nature, and we therefore need that the perturbation terms $\mathcal{O}(|n-1|/r)$ appearing in (\ref{SHEqn2}) are small. Thus, our results focus on the case $n=1+\varepsilon$ with $0<\varepsilon\ll1$ (note that we can consider $n$ as a real parameter in (\ref{SHEqn2}) though $n$ is then no longer related to the space dimension) and on the case $n=2,3$ with $r\gg1$ large. We now outline our results:
\begin{compactitem}
\item For $|n-1|\ll1$, we show that snaking branches persist for plateau lengths $L\leq\exp(b/|n-1|)$ where $b>0$ is a constant (Theorem~\ref{thm:epsilonsmall}).
\item For $|n-1|\ll1$, we will study under which conditions on the perturbation terms localized rolls cannot persist for large plateau lengths $L\geq L_\mathrm{max}(|n-1|)$ and when they will persist for all large $L\gg1$ (Theorem~\ref{thm:collapsedsnaking}).
\item For $n=2,3$, we will give conditions on the perturbation terms under which localized rolls cannot persist for large plateau lengths $L\gg1$ (Theorem~\ref{thm:epsilonlarge}).
\item For the planar and three-dimensional Swift--Hohenberg equation, we will show using analytical and numerical results that snaking branches need to collapse onto the Maxwell point (\S\ref{sec:SwiftHohenberg}).
\end{compactitem}
We emphasize that our results will be formulated for a general class of systems that includes (\ref{SHEqn2}).

The remainder of this paper is organized as follows. We summarize our hypotheses and main results in \S\ref{sec:Results} and apply these results to the Swift--Hohenberg equation in \S\ref{sec:SwiftHohenberg}. The remaining sections are dedicated to the proofs of our main theorems. We will construct boundary-layer solutions near the singularity of (\ref{SHEqn2}) at $r=0$ in \S\ref{sec:BoundaryLayer}, discuss the dynamics near the family of periodic orbits in \S\ref{sec:Fenichel}, consider the stable manifold of $u=0$ in \S\ref{sec:StableManifold}, and construct radial pulses in \S\ref{sec:Matching}. In \S\ref{sec:IntegralManifold}, we expand the vector field on the integral manifold and use these results in \S\ref{sec:persistent} to analyse when snaking persists and when collapsed snaking occurs.


\section{Main results} \label{sec:Results}

Consider the ordinary differential equation
\begin{equation} \label{ODE}
u_x = f(u,\mu),
\end{equation}
where $u \in \mathbb{R}^4$, $\mu \in \mathbb{R}$, and $f:\mathbb{R}^4 \times \mathbb{R} \to \mathbb{R}^4$ is smooth. Our first assumption concerns reversibility.

\begin{hypothesis} \label{hyp:Reverser}
There exists a linear map $\mathcal{R}:\mathbb{R}^4 \to \mathbb{R}^4$ with $\mathcal{R}^2 = 1$ and $\mathrm{dim\ Fix}(\mathcal{R}) = 2$ so that $f(\mathcal{R}u,\mu) = -\mathcal{R}f(u,\mu)$ for all $(u,\mu)$.
\end{hypothesis}

Hypothesis~\ref{hyp:Reverser} implies that if $u(x)$ is a solution to (\ref{ODE}), then so is $\mathcal{R}u(-x)$. Furthermore, if $u(0)\in\mathrm{Fix}(\mathcal{R})$ we have that $u(x)=\mathcal{R}u(-x)$ for all $x \in \mathbb{R}$, and hence we refer to such solutions as \emph{symmetric}. Finally, we remark that $\mathbb{R}^4 = \mathrm{Fix}(\mathcal{R}) \oplus \mathrm{Fix}(-\mathcal{R})$. Next, we assume the existence of a conserved quantity.

\begin{hypothesis} \label{hyp:Energy}
There exists a smooth function $\mathcal{H}:\mathbb{R}^4\times\mathbb{R} \to \mathbb{R}$ with $\mathcal{H}(\mathcal{R}u,\mu) = \mathcal{H}(u,\mu)$ and $\langle \nabla_u \mathcal{H}(u,\mu),f(u,\mu)\rangle = 0$ for all $(u,\mu)$. We normalize $\mathcal{H}$ so that $\mathcal{H}(0,\mu) = 0$ for all $\mu$.
\end{hypothesis}

Our next hypothesis states that the origin is a hyperbolic saddle.

\begin{hypothesis} \label{hyp:TrivialEquib}
We assume that $f(0,\mu) = 0$ for all $\mu$ and that $f_u(0,\mu)$ has exactly two eigenvalues with strictly negative real part and two eigenvalues with strictly positive real part.
\end{hypothesis}

Next, we formalize the existence of hyperbolic periodic orbits that are parametrized by the value of the conserved quantity $\mathcal{H}(\cdot,\mu)$. Throughout this paper, we denote the interior of an interval $J$ by $\mathring{J}$.

\begin{hypothesis} \label{hyp:Periodic}
There exist compact intervals $J,K\subset\mathbb{R}$ with $\mathring{J}\neq\emptyset$ and $0\in\mathring{K}$ such that (\ref{ODE}) has, for each $(\mu,h)\in J\times K$, a periodic orbit $\gamma(x,\mu,h)$ with minimal period $p(\mu,h) > 0$ such that the following holds for each $(\mu,h)\in J\times K$:
\begin{compactenum}
\item $\gamma(x,\mu,h)$ and $p(\mu,h)$ depend smoothly on $(\mu,h)$.
\item $\gamma(x,\mu,h)$ is symmetric: $\gamma(0,\mu,h) \in \mathrm{Fix}(\mathcal{R})$.
\item $\mathcal{H}(\gamma(x,\mu,h),\mu)=h$ and $\mathcal{H}_u(\gamma(x,\mu,h),\mu) \neq 0$ for one, and hence all, $x$.
\item Each $\gamma(x,\mu,h)$ has two positive Floquet multipliers $e^{\pm\alpha(\mu,h)p(\mu,h)}$ that depend smoothly on $(\mu,h)$ and satisfy $\min_{(\mu,h)\in J\times K} \alpha(\mu,h)>0$.
\end{compactenum}
\end{hypothesis}

Reversibility implies that the set of Floquet exponents of a symmetric periodic orbit is invariant under multiplication by $-1$. As shown in \cite{Aougab}, the case where the two hyperbolic Floquet multipliers are negative may not lead to snaking. Hypothesis~\ref{hyp:Periodic} implies that the union $\mathcal{P}(\mu):=\{\gamma(x,\mu,h):\, x\in\mathbb{R}, h\in K\}$ of the periodic orbits is a normally hyperbolic invariant manifold parametrized by the value $h\in K$ of the conserved quantity.

As in \cite{Beck}, we restrict the system (\ref{ODE}) to the three-dimensional level set $\mathcal{H}^{-1}(0)$ and parametrize a neighborhood of the periodic orbit $\gamma(\cdot,\mu,0)$ using the variables $(\varphi,v^s,v^u)$, where $(\varphi,0,0)$ corresponds to $\gamma(\varphi p(\mu,0)/2\pi,\mu,0)$, and $(\varphi,v^s,0)$ and $(\varphi,0,v^u)$ parametrize the strong stable and strong unstable fibers $W^{ss}(\gamma(\varphi p(\mu,0)/2\pi,\mu,0),\mu)$ and $W^{uu}(\gamma(\varphi p(\mu,0)/2\pi,\mu,0),\mu)$, respectively, of $\gamma(\varphi p(\mu,0)/2\pi,\mu,0)$. Using the coordinates $(\varphi,v^s,v^u)$, we then define the section
\[
\Sigma_\mathrm{out} := \{ (\varphi,v^s,v^u,h)\in S^1\times[-\delta,\delta]\times[-\delta,\delta]\times K:\, v^u=\delta \},
\]
where $\delta>0$ is a small positive constant. We can now formulate our assumptions on the existence of heteroclinic orbits that connect the periodic orbits $\gamma$ to the rest state $u=0$.

\begin{hypothesis} \label{hyp:StableIntersection}
There exists a smooth function $G:S^1\times\mathcal{I}\times J\to\mathbb{R}$ such that $G(\varphi,v^s,\mu)=0$ if and only if $(\varphi,v^s,\delta,0)\in W^s(0,\mu)\cap\Sigma_\mathrm{out}$. In particular,
\begin{eqnarray} \label{Gamma}
\Gamma & := & \{(\varphi,\mu)\in S^1\times J:\, G(\varphi,0,\mu) = 0\} \\ \nonumber & = &
\{(\varphi,\mu)\in S^1\times J:\, W^s(0,\mu)\cap W^{uu}(\gamma(\varphi p(\mu,0)/2\pi,\mu,0),\mu)\cap\Sigma_\mathrm{out}\neq\emptyset \},
\end{eqnarray}
and we assume that $\Gamma\subset S^1\times\mathring{J}$ is nonempty with $D_{(\varphi,\mu)} G(\varphi,0,\mu)\neq0$ for each $(\varphi,\mu)\in\Gamma$.
\end{hypothesis}

\begin{figure}
\centering
\includegraphics{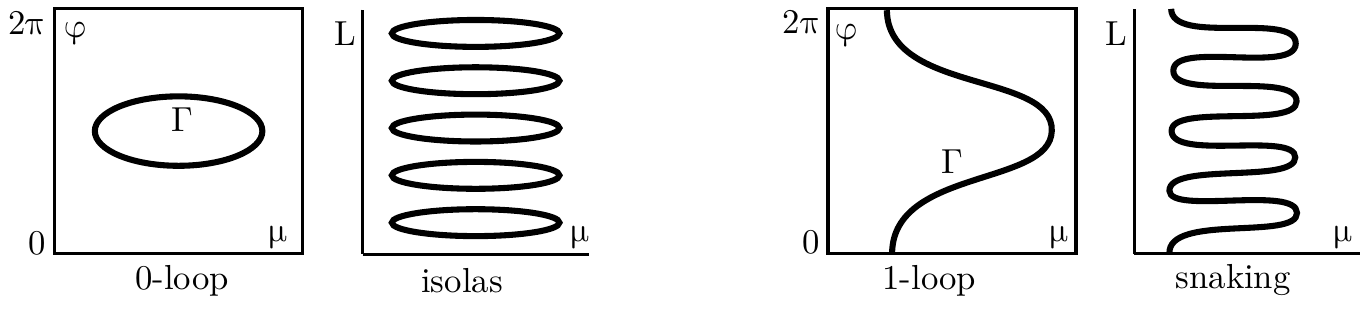}
\caption{Shown are different possible configurations of the set $\Gamma$: The two panels on the left show a 0-loop and the resulting isolas. In contrast, as shown in the two rightmost panels, 1-loops lead to snaking diagrams.}
\label{f:5}
\end{figure}

As shown in \cite{Beck,Aougab}, Hypothesis~\ref{hyp:StableIntersection} implies that $\Gamma$ is the union of finitely many disjoint closed loops. Parametrizing one such loop by a function $(\varphi(s),\mu(s))$ with $s\in[0,1]$ and $\varphi(s)$ in the universal cover $\mathbb{R}$ of $S^1$, we have either (i) $\varphi(0)=\varphi(1)$ or (ii) $\varphi(0)\neq\varphi(1)$. Following \cite{Aougab}, we will refer to the case (i) as a 0-loop and case (ii) as a 1-loop. As proved in \cite{Beck,Aougab} and illustrated in Figure~\ref{f:5}, 0-loops lead to isolas and 1-loops to snaking branches. We denote by $\Gamma_\mathrm{lift}\subset\mathbb{R}\times\mathring{J}$ the preimage of $\Gamma$ under the natural covering projection from $\mathbb{R}\times\mathring{J}$ to $S^1\times\mathring{J}$ so that 0-loops in $\Gamma$ are lifted to an infinite number of disjoint copies of the 0-loop, whereas 1-loops lift to an unbounded connected curve.

Motivated by the structure of (\ref{SHEqn2}), our goal is to extend the results in \cite{Beck} to systems of the form
\begin{equation} \label{ODEPert}
u_x = f(u,\mu) + \frac{\varepsilon}{x}g(u,\mu,\varepsilon),
\end{equation}
where $\varepsilon>0$ is not necessarily small.

\begin{hypothesis} \label{hyp:Pert}
The function $g:\mathbb{R}^4\times\mathbb{R}\times\mathbb{R}^+\to \mathbb{R}^4$ is smooth in all its arguments, and $g(u,\mu,\varepsilon)=0$ for all $(u,\mu,\varepsilon)\in\mathrm{Fix}(\mathcal{R})\times J\times\mathbb{R}^+$.
\end{hypothesis}

Hypothesis~\ref{hyp:Pert} implies in particular that $u=0$ is a solution of (\ref{ODEPert}) for all values of $\varepsilon$. We are interested in constructing solutions to $(\ref{ODEPert})$ that remain close to the manifold $\mathcal{P}(\mu)$ of periodic orbits for $x\in[0,L]$ for appropriate large values of $L\gg1$ and converge to $u=0$ as $x\to\infty$. To make this more precise, we denote by $U_\delta(\mathcal{P}(\mu))$ the $\delta$-neighbourhood of the manifold $\mathcal{P}(\mu)$ and by $W_L^s(0,\mu,\varepsilon)\subset\mathbb{R}^4$ the slice of the stable manifold of the rest state $u=0$ of (\ref{ODEPert}) for $x=L$. We then say that $u(x)$ is a \emph{radial pulse} with plateau length $L$ for some $L\gg1$ if $u(x)$ is defined for $x\geq0$, is a solution of (\ref{ODEPert}) for $x>0$ with $(\mu,\varepsilon)$ fixed, and satisfies the conditions
\begin{equation}\label{bc1}
u(0)\in\mathrm{Fix}(\mathcal{R}), \qquad
u(x)\in U_\delta(\mathcal{P}(\mu)) \mbox{ for } x\in[0,L],\qquad
u(L)\in \partial U_\delta(\mathcal{P}(\mu))\cap W_L^s(0,\mu,\varepsilon);
\end{equation}
see Figure~\ref{fig:Cylinder} for an illustration. Our first result relates the structure of $\Gamma_\mathrm{lift}$ to the bifurcation structure of radial pulses when $0<\varepsilon\ll1$.

\begin{theorem} \label{thm:epsilonsmall}
Assume that Hypotheses~\ref{hyp:Reverser}-\ref{hyp:Pert} are met, then there are constants $b,\varepsilon_0,\eta,L_*>0$, a function $L_\mathrm{max}(\varepsilon)$ such that $L_\mathrm{max}(\varepsilon)\geq \mathrm{e}^{b/\varepsilon}$, and sets $\Gamma_\mathrm{pulse}^{\varphi_0,\varepsilon}\subset (L_*,L_\mathrm{max}(\varepsilon)) \times J$ defined for $\varphi_0\in\{0,\pi\}$ and $\varepsilon\in[0,\varepsilon_0]$ so that the following is true:
\begin{compactenum}
\item Equation (\ref{ODEPert}) admits a radial pulse if and only if $(L,\mu)\in\Gamma_\mathrm{pulse}^{\varphi_0,\varepsilon}$ for $\varphi_0=0$ or $\varphi_0=\pi$.
\item There exists a smooth function $\tilde{g}_\mathrm{lift}(L,\mu,\varepsilon) = \mathcal{O}(\varepsilon \ln L)$ such that for each fixed $\varphi_0\in\{0,\pi\}$ and $\varepsilon\in[0,\varepsilon_0]$ the one-dimensional manifolds
\[
\tilde{\Gamma}_\mathrm{lift}^{\varphi_0,\varepsilon} := \{(L-\tilde{g}_\mathrm{lift}(L,\mu,\varepsilon)+\varphi_0,\mu) :\ (L,\mu) \in \Gamma_\mathrm{lift}\cap((L_*,L_\mathrm{max}(\varepsilon)) \times J)\}
\]
and $\Gamma_\mathrm{pulse}^{\varphi_0,\varepsilon}$ are $\mathcal{O}(\mathrm{e}^{-\eta L})$-close to each other in the $C^0$-sense near each point $(L,\mu)\in\tilde{\Gamma}_\mathrm{lift}^{\varphi_0,\varepsilon}$.
\end{compactenum}
\end{theorem}

We emphasize that Theorem~\ref{thm:epsilonsmall} captures not only those solutions that stay close to the level set $\mathcal{H}^{-1}(0)$ but also all solutions along which the function $\mathcal{H}$ takes values in the interval $K$. In particular, the size of $b$ is restricted only by the possibility that a solution leaves a neighborhood of the manifolds $\mathcal{P}(\mu)$ when the value of the quantity $\mathcal{H}$ reaches the boundary of the interval $K$.

Our next result gives conditions for collapsed snaking for $0\leq\varepsilon\ll1$. To state the theorem, we define the function
\begin{equation}\label{e:s}
S(h,\mu) := \frac{1}{p(\mu,h)}\int_0^{p(\mu,h)} \langle \nabla_u \mathcal{H}(\gamma(x,\mu,h)),g(\gamma(x,\mu,h),\mu,0)\rangle\,\mathrm{d}x,
\end{equation}
which is equal to the average of the perturbation $g$ in the direction of the gradient of $\mathcal{H}$ along the periodic orbits. We will see in \S\ref{sec:IntegralManifold} that $S(h,\mu)$ is the vector field that describes, to leading order, via the differential equation
\begin{equation}\label{e:hvf}
h_x = \frac{\varepsilon}{x} S(h,\mu)
\end{equation}
how the value $h(x)$ of the conserved quantity $\mathcal{H}(u,\mu)$ changes along solutions $u(x)$ of (\ref{ODEPert}). A necessary condition for the existence of radial pulses $u(x)$ with plateau length $L$ is that $h(x)\in K$ for $0\leq x\leq L$ and $h(L)\approx0$, as the latter is necessary for $u(x)$ to satisfy $u(L)\in W^s_L(0,\mu)$. Our next theorem states conditions on the vector field $S(h,\mu)$ that preclude or guarantee that solutions $h(x)$ of (\ref{e:hvf}) stay in $K$ for all $x\geq0$.

\begin{theorem} \label{thm:collapsedsnaking}
Assume that Hypotheses~\ref{hyp:Reverser}-\ref{hyp:Pert} are met.
\begin{compactenum}
\item If there is a closed interval $\tilde{J}\subset J$ such that $S(h,\mu)<0$ for all $h\in K\cap\mathbb{R}^+$ and $\mu\in\tilde{J}$ (or, alternatively, $S(h,\mu)>0$ for all $h\in K\cap\mathbb{R}^-$ and $\mu\in\tilde{J}$), then there are a constant $\varepsilon_0>0$ and a function $L_\mathrm{min}(\varepsilon)$ so that (\ref{ODEPert}) with $\mu\in\tilde{J}$ and $0<\varepsilon<\varepsilon_0$ cannot have any radial pulses with plateau lengths $L\geq L_\mathrm{min}(\varepsilon)$.
\item Assume that there are $\varphi\in S^1$ and $\mu_*\in\mathring{J}$ such that $S(0,\mu_*)=0$, $S_h(0,\mu_*)<0$, $S_\mu(0,\mu_*)>0$, $(\varphi,\mu_*)\in\Gamma$, and $G_\varphi(\varphi,0,\mu_*)\neq0$, then there exists an $\varepsilon_0>0$ such that the following is true for each $\varphi_0\in\{0,\pi\}$ and each $0<\varepsilon<\varepsilon_0$: there exists a sequence $(L_m,\mu_m)$ with $L_m\to\infty$ monotonically as $m\to\infty$ and $\mu_m$ near $\mu_*$ for all $m$ so that (\ref{ODEPert}) with $\mu=\mu_m$ has a radial pulse with plateau length $L_m$.
\end{compactenum}
\end{theorem}

Next, we focus on arbitrary, not necessarily small values of $\varepsilon>0$. We say that $u(x)$ is an \emph{$R$-asymptotic radial pulse} of plateau length $L$ if $u(x)$ is defined for $x\geq R$, is a solution of (\ref{ODEPert}) for $x>R$ with $(\mu,\varepsilon)$ fixed, and satisfies the conditions
\[
u(x)\in U_\delta(\mathcal{P}(\mu)) \mbox{ for } x\in[R,L],\qquad
u(L)\in \partial U_\delta(\mathcal{P}(\mu))\cap W_L^s(0,\mu,\varepsilon).
\]
Implicit in our definition is the assumption that $L>R$. Our next result provides conditions on the existence and nonexistence of $R$-asymptotic radial pulses. In contrast to our definition of radial pulses in (\ref{bc1}), we do not impose any boundary conditions for $R$-asymptotic radial pulses at $x=0$ or $x=R$ and can therefore guarantee the existence of these solutions for all sufficiently large $L$ instead of just for a sequence as in Theorem~\ref{thm:collapsedsnaking}.

\begin{theorem} \label{thm:epsilonlarge}
Assume that Hypotheses~\ref{hyp:Reverser}-\ref{hyp:Pert} are met.
\begin{compactenum}
\item If there is a closed interval $\tilde{J}\subset J$ such that $S(h,\mu)<0$ for all $h\in K\cap\mathbb{R}^+$ and $\mu\in\tilde{J}$ (or, alternatively, $S(h,\mu)>0$ for all $h\in K\cap\mathbb{R}^-$ and $\mu\in\tilde{J}$), then for each fixed and not necessarily small $\varepsilon_*>0$ there are constants $1<R_*<L_*<\infty$ such that (\ref{ODEPert}) with $\mu\in\tilde{J}$ and $\varepsilon=\varepsilon_*$ cannot have any $R_*$-asymptotic radial pulses with plateau lengths $L\geq L_*$.
\item If there are constants $\varphi\in S^1$ and $\mu_*\in\mathring{J}$ such that $S(0,\mu_*)=0$, $S_h(0,\mu_*)<0$, $S_\mu(0,\mu_*)>0$, $(\varphi,\mu_*)\in\Gamma$, and $G_\varphi(\varphi,0,\mu_*)\neq0$, then for each fixed and not necessarily small $\varepsilon_*>0$ and each $\delta>0$ there are constants $1<R_*<L_*<\infty$ and a function $\mu(L)$ defined for $L\geq L_*$ with $|\mu(L)-\mu_*|<\delta$ such that (\ref{ODEPert}) with $\varepsilon=\varepsilon_*$ has an $R_*$-asymptotic radial pulse with plateau length $L$ at $\mu=\mu(L)$ for each $L\geq L_*$.
\end{compactenum}
\end{theorem}

We will see in \S\ref{sec:SwiftHohenberg} that if $\mu_*$ is the Maxwell point and $\tilde{J}$ is any closed interval in $J\setminus\{\mu_*\}$, then the Swift--Hohenberg equation satisfies the conditions stated in Theorem~\ref{thm:epsilonlarge}, and snaking therefore has to collapse onto the Maxwell point for $n=2,3$. Theorem~\ref{thm:epsilonsmall} will be proved in \S\ref{sec:Matching}, while Theorems~\ref{thm:collapsedsnaking} and~\ref{thm:epsilonlarge} will be proved in \S\ref{sec:persistent}.


\section{Application to the Swift--Hohenberg equation} \label{sec:SwiftHohenberg}

We now apply the results presented in the preceding section to the Swift--Hohenberg equation
\begin{equation}\label{e:SHn}
U_t = -(1+\Delta)^2 U - \mu U + \nu U^2 - U^3, \qquad x\in\mathbb{R}^n.
\end{equation}
Radial solutions of this equation satisfy the PDE
\begin{equation} \label{SHEqn}
0 = -\left(1 + \frac{n-1}{r}\partial_r + \partial_r^2\right)^2U - \mu U + \nu U^2 - U^3,
\end{equation}
where $r=|x|$ denotes the radial direction in $\mathbb{R}^n$. Throughout this section, we will keep $\nu$ fixed and vary $\mu$: in particular, we will not explicitly indicate the dependence of any quantities on $\nu$.

Using a combination of analytical and numerical results, we will show that the Swift--Hohenberg equation satisfies the assumptions stated in \S\ref{sec:Results} and that the snaking branches for $n>1$ have to collapse onto a single value $\mu_\mathrm{Max}$ of the parameter $\mu$ as $L\to\infty$. We also identify $\mu_\mathrm{Max}$ with the Maxwell point.

\subsection{Verification of Hypotheses~\ref{hyp:Reverser}--\ref{hyp:Pert}}

\begin{figure}
\centering
\includegraphics{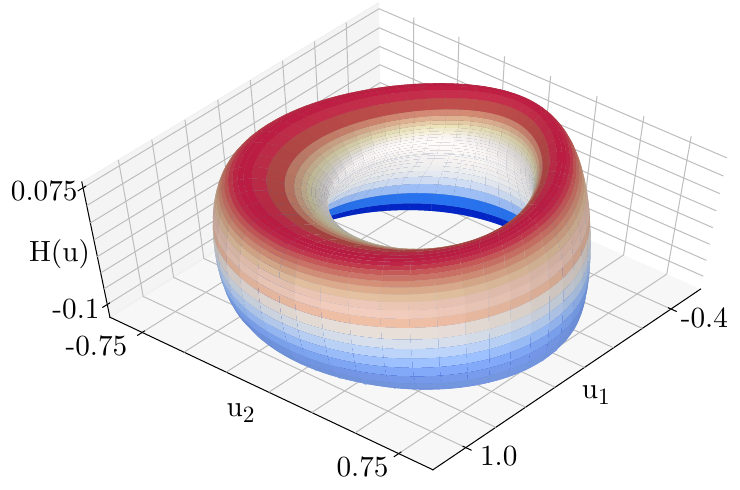}
\caption{Shown is the torus of periodic orbits of (\ref{FirstOrderSH}) for $\mu=0.2004$, $\nu=1.6$, and $\varepsilon=0$. The solutions on the inside of the torus are hyperbolic, while those on the outside are elliptic.}
\label{fig:Torus}
\end{figure}

We define $\varepsilon:=n-1$, $u_1=U$, $u_2=\partial_r u_1$, $u_3=(1+\frac{\varepsilon}{r}\partial_r + \partial_r^2)u_1$, and $u_4=\partial_r u_3$, then (\ref{SHEqn}) can be written as the first-order system
\begin{equation} \label{FirstOrderSH}
\begin{split}
u_1^\prime &= u_2, \\
u_2^\prime &= u_3 - u_1 - \frac{\varepsilon}{r} u_2, \\
u_3^\prime &= u_4, \\
u_4^\prime &= -u_3 -\mu u_1 + \nu u_1^2 - u_1^3 - \frac{\varepsilon}{r} u_4,
\end{split}
\end{equation}
where $\prime$ denotes differentiation with respect to $r$. Setting $\varepsilon=0$ in (\ref{FirstOrderSH}), we find that the resulting system is reversible with reverser
\[
\mathcal{R} = \begin{pmatrix} 1 & 0 & 0 & 0 \\ 0 & -1 & 0 & 0 \\
0 & 0 & 1 & 0 \\ 0 & 0 & 0 & -1 \end{pmatrix}
\]
and has the conserved quantity
\begin{equation}\label{EnergyH}
\mathcal{H}(u,\mu) = u_2 u_4 + u_1 u_3 -\frac{u_3^2}{2} + \frac{\mu u_1^2}{2}
- \frac{\nu u_1^3}{3} + \frac{u_1^4}{4}.
\end{equation}
It is now straightforward to verify that Hypotheses~\ref{hyp:Reverser}--\ref{hyp:TrivialEquib} hold. Figure~\ref{fig:Torus} reflects the numerical evidence for the existence of a torus of periodic orbits $\gamma(r,\mu,h)$ to (\ref{FirstOrderSH}) when $\varepsilon=0$. Numerically, the periodic orbits in the inside of the torus shown in Figure~\ref{fig:Torus} are hyperbolic, thus indicating that Hypothesis~\ref{hyp:Periodic} is indeed met. Furthermore, Figure~\ref{f:3}(i) contains the numerical snaking diagram of localized rolls of (\ref{FirstOrderSH}) for $\varepsilon=0$. As shown in \cite{Beck}, the structure of the branches visible in this figure is consistent with the assumption that the set $\Gamma$ consists of a single 1-loop that satisfies Hypothesis~\ref{hyp:StableIntersection}. Finally, allowing $\varepsilon\neq0$, we see that
\[
g(u,\mu,\varepsilon) = -\begin{pmatrix} 0 \\ u_2 \\ 0 \\ u_4 \end{pmatrix},
\]
which vanishes precisely when $u\in\mathrm{Fix}(\mathcal{R})$ as required in Hypothesis~\ref{hyp:Pert}. With the caveat that Hypotheses~\ref{hyp:Periodic}--\ref{hyp:StableIntersection} can be verified only numerically, Theorem~\ref{thm:epsilonsmall} implies that snaking persists for each fixed $0<\varepsilon\ll1$ for all $L$ with $L\leq\mathrm{e}^{b/\varepsilon}$ for some constant $b>0$.

\subsection{Collapsed snaking}

Next, we use a combination of analytical and numerical results to show that the snaking branches for the Swift--Hohenberg equation (\ref{e:SHn}) with $n>1$ need to collapse onto a single value of $\mu$ as $L\to\infty$. The one-dimensional Swift--Hohenberg equation
\begin{equation}\label{e:sh}
U_t = -(1+\partial_x^2)^2 U - \mu U + \nu U^2 - U^3, \qquad x\in\mathbb{R}
\end{equation}
considered on the space of $p$-periodic functions admits the PDE energy functional
\begin{equation}\label{Energy}
\mathcal{E}(U,\mu,p) := \frac{1}{p} \int_0^p \left( \frac{(U+U_{xx})^2}{2} + \frac{\mu U^2}{2} - \frac{\nu U^3}{3} + \frac{U^4}{4} \right) \,\mathrm{d}x.
\end{equation}
Our goal is to relate the function $S(h,\mu)$ defined in (\ref{e:s}) to the energy functional $\mathcal{E}$ and the conserved quantity $\mathcal{H}$. Before stating our result, we introduce additional notation. We denote by $U^*(x,\mu,h):=\gamma_1(x,\mu,h)$ the stationary roll solutions of (\ref{e:sh}) with minimal spatial period $p(\mu,h)>0$ that satisfy $\mathcal{H}(U^*(x,\mu,h),\mu)=h$ for one, and hence all, $x$. We say that $\mu=\mu_\mathrm{Max}$ is a Maxwell point if $\mathcal{E}(U^*(\cdot,\mu,0),\mu,p(\mu,0))=0$, so that the PDE energy of the roll solution with $\mathcal{H}=0$ vanishes. Numerically, (\ref{e:sh}) has a unique Maxwell point $\mu=\mu_\mathrm{Max}$ for each value of $\nu$.

We can now calculate the function $S(h,\mu)$ defined in (\ref{e:s}), which, via the differential equation
\begin{equation}\label{e:hvf1}
h_x = \frac{\varepsilon}{x} S(h,\mu),
\end{equation}
describes to leading order how the value $h(x)$ of $\mathcal{H}(U(x),\mu)$ changes along a radial pulse $U(x)$. As pointed out in \S\ref{sec:Results}, a necessary condition for the existence of radial pulses is that $h(x)\in K$ for all $0\leq x\leq L$ and $h(L)\approx0$. Using (\ref{e:s}) and the form of $\mathcal{H}$ and $g$ discussed in the last section, we find
\begin{eqnarray*}
S(h,\mu) & = &
\frac{1}{p(\mu,h)}\int_0^{p(\mu,h)} \langle \nabla_u \mathcal{H}(\gamma(x,\mu,h)),g(\gamma(x,\mu,h),\mu,0)\rangle \,\mathrm{d}x \\ & = &
\frac{-2}{p(\mu,h)}\int_0^{p(\mu,h)} \gamma_2(x,\mu,h) \gamma_4(x,\mu,h) \,\mathrm{d}x \\ & = &
\frac{-2}{p(\mu,h)}\int_0^{p(\mu,h)} (U^*_x(x,\mu,h)^2 + U^*_x(x,\mu,h) U^*_{xxx}(x,\mu,h)) \,\mathrm{d}x.
\end{eqnarray*}
Our main result relates the vector field $S(h,\mu)$ to the energy $\mathcal{E}$.

\begin{lemma}\label{lem:PDE_Energy}
We have
\[
S(h,\mu) = \mathcal{E}(U^*(\cdot,\mu,h),\mu,p(\mu,h)) - h, \qquad
S_h(0,\mu_\mathrm{Max}) < 0, \qquad
S_\mu(0,\mu_\mathrm{Max}) > 0.
\]
In particular, $S(0,\mu)=0$ if and only if $\mu=\mu_\mathrm{Max}$.
\end{lemma}

\begin{figure}
\centering
\includegraphics[scale=1]{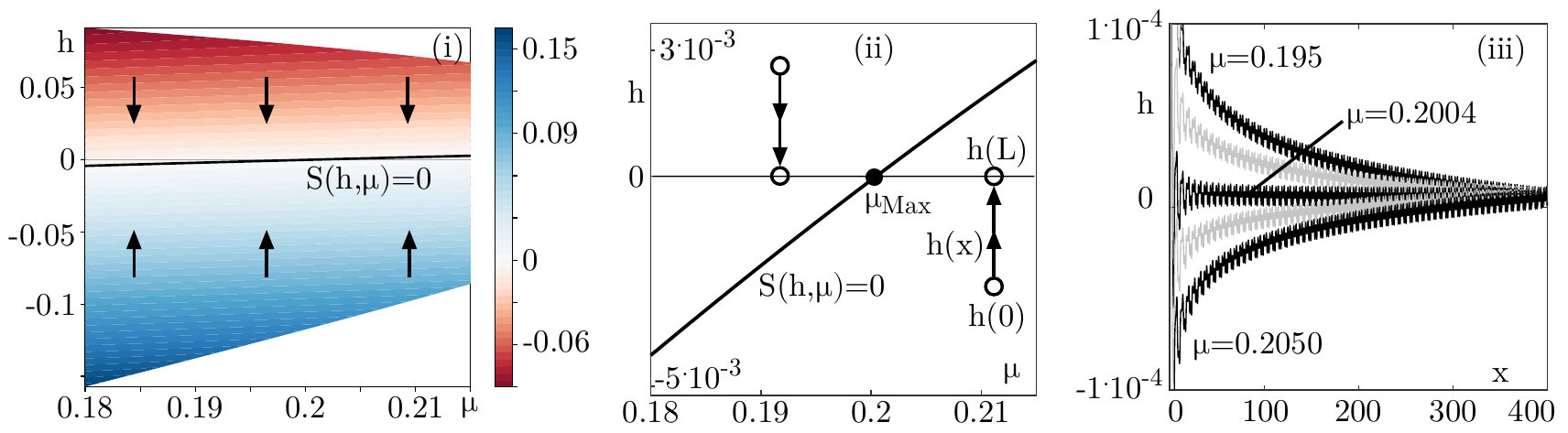}
\caption{Shown are \textsc{auto09p} computations for the Swift--Hohenberg equation (\ref{e:sh}) with $\nu=1.6$. Panel~(i) shows a contour plot of the function $S(h,\mu)$, with arrows indicating the directions in which solutions $h(x)$ of (\ref{e:hvf1}) would move. An enlarged plot of the curve $S(h,\mu)=0$ is shown in panel~(ii), where we also include two schematic trajectories of solutions $h(x)$ of (\ref{e:hvf1}) for $0\leq x\leq L$ that satisfy $h(L)=0$ with $L\gg1$. Panel~(iii) contains the graphs of $\mathcal{H}(U(x),\mu)$ along radial pulses at $\varepsilon=0.1$ for $\mu=0.195,0.1975,0.2004,0.2025,0.205$ (corresponding to the curves from top to bottom), corroborating the theoretical prediction for solutions $h(x)=\mathcal{H}(U(x),\mu)$ of (\ref{e:hvf1}) depending on whether $\mu$ is smaller or larger than the Maxwell point $\mu_\mathrm{Max}=0.2004$.}
\label{fig:AvgEnergy}
\end{figure}

Before proving this result, we discuss its implications for the Swift--Hohenberg equation (\ref{e:SHn}) posed on $\mathbb{R}^n$.

First, Lemma~\ref{lem:PDE_Energy} shows that the Swift--Hohenberg equation satisfies the hypotheses on $S$ needed in  Theorem~\ref{thm:collapsedsnaking}(ii), while Figure~\ref{f:3}(i) indicates that the assumptions on $\Gamma$ are met. For $|n-1|\ll1$, we can therefore conclude that radial pulses with arbitrarily large plateau lengths exist for parameter values $\mu$ near the Maxwell point.

Next, we can use \textsc{auto09p} to compute the vector field $S(h,\mu)$ numerically through continuation of the periodic solutions of the one-dimensional Swift--Hohenberg equation. The results shown in Figure~\ref{fig:AvgEnergy}(i) and~(ii) indicate that (\ref{e:hvf1}) has a unique equilibrium for each value of $\mu$, and that these equilibria are stable---note that Lemma~\ref{lem:PDE_Energy} provides a proof of these properties, including the location of the equilibria, for $\mu$ near the Maxwell point. In particular, these numerical results obtained for $n=1$ show that the hypotheses of Theorem~\ref{thm:epsilonlarge}(i) are met for each $n>1$ (including $n=2,3$), and we conclude that $R$-asymptotic radial pulses of plateau length $L$ cannot exist for $L\gg1$ and that the snaking branches therefore have to collapse onto the Maxwell point for each $n>1$. We argued above that the assumptions for Theorem~\ref{thm:epsilonlarge}(ii) are also met, and we can conclude that $R_*$-asymptotic radial pulses exist near the Maxwell point for arbitrarily large plateau lengths $L$. Note that our definition of $R_*$-asymptotic radial pulses ignores the spatial interval $[0,R_*]$: our results therefore apply equally to branches involving spots and rings, but they cannot make any predictions for actual radial pulses as this would require that we construct solutions on $[0,R_*]$ and match them with the $R_*$-asymptotic pulses at $x=R_*$.

Finally, as illustrated in the schematic in Figure~\ref{fig:AvgEnergy}(ii), the solutions $u(x)$ that reach the $\mathcal{H}=0$ level set at $x=L$ necessarily have $\mathcal{H}(u(0),\mu)>0$ for $\mu<\mu_\mathrm{Max}$ and $\mathcal{H}(u(0),\mu)<0$ for $\mu>\mu_\mathrm{Max}$, provided $L$ is sufficiently large. Figure~\ref{fig:AvgEnergy}(iii) confirms this prediction. We now give the proof of Lemma~\ref{lem:PDE_Energy}.

\begin{proof}[Proof of Lemma~\ref{lem:PDE_Energy}]
To prove the characterization of $S(h,\mu)$ in terms of the energy $\mathcal{E}$, we use the notation $U(x):=U^*(x,\mu,h)$ and $p:=p(\mu,h)$. Note that $U(x)$ is then $p$-periodic. Writing the conserved quantity $\mathcal{H}$ defined in (\ref{EnergyH}) in terms of derivatives of $U$, we find that
\[
\mathcal{H}(U,\mu) = U_x U_{xxx} - \frac{U^2_{xx}}{2} + U^2_x + \frac{U^2}{2} + \frac{\mu U^2}{2} - \frac{\nu U^3}{3} + \frac{U^4}{4}
\]
is conserved pointwise along $U(x)$. Writing the pointwise identity $\mathcal{H}(U(x),\mu)=h$ as
\[
\frac{\mu U^2}{2} - \frac{\nu U^3}{3} + \frac{U^4}{4} = h - U_x U_{xxx} + \frac{U^2_{xx}}{2} - U^2_x - \frac{U^2}{2},
\]
substituting this identity into (\ref{Energy}), and integrating by parts gives
\begin{eqnarray*}
\mathcal{E}(U(\cdot),\mu,p) & = &
\frac{1}{p} \int_0^p \left( \frac{(U+U_{xx})^2}{2} + \frac{\mu U^2}{2} - \frac{\nu U^3}{3} + \frac{U^4}{4} \right) \mathrm{d}x \\ & = &
\frac{1}{p} \int_0^p \left( \frac{U^2}{2} + U U_{xx} + \frac{U^2_{xx}}{2} + h - U_x U_{xxx} + \frac{U^2_{xx}}{2} - U^2_x - \frac{U^2}{2} \right) \mathrm{d}x \\ & = &
h + \frac{1}{p} \int_0^p \left( U U_{xx} + U^2_{xx} - U_x U_{xxx} - U^2_x \right)\mathrm{d}x \\ & = &
h - \frac{2}{p} \int_0^p \left(U_x^2 + U_x U_{xxx}\right)\mathrm{d}x \\ & = &
h + S(h,\mu)
\end{eqnarray*}
as claimed. Next, we prove the claims about the derivatives. First, we consider the expression for the energy and rescale $x=py$ to get
\[
\mathcal{E}(V(\cdot),\mu,p) = \int_0^1 \left( \frac12 \left(V+\frac{V_{yy}}{p^2}\right)^2 + \frac{\mu V^2}{2} - \frac{\nu V^3}{3} + \frac{V^4}{4} \right) \mathrm{d}y,
\]
where $V(y)=U(py)$ is 1-periodic in $y$. Taking the derivative of this expression with respect to $h$ and using that $\mathcal{E}_V(V(\cdot),\mu,p)=0$, we obtain
\begin{eqnarray*}
\frac{\mathrm{d}}{\mathrm{d}h} \left( \mathcal{E}(V(\cdot),\mu,p) \right) & = &
\mathcal{E}_V(V(\cdot),\mu,p) V_h(\cdot) - \frac{2p_h}{p^3} \int_0^1 \left(V+\frac{V_{yy}}{p^2}\right) V_{yy}\, \mathrm{d}y \;=\;
- \frac{2p_h}{p^2} \int_0^p (U+U_{xx}) U_{xx}\, \mathrm{d}x \\ & = &
\frac{2p_h}{p^2} \int_0^p (U_x+U_{xxx}) U_x\, \mathrm{d}x \;=\; 
- \frac{p_h}{p} S(h,\mu)
\end{eqnarray*}
and therefore
\[
S_h(0,\mu_\mathrm{Max})
= \frac{\mathrm{d}}{\mathrm{d}h} \left( \mathcal{E}(V(\cdot),\mu,p) \right)\Big|_{(h,\mu)=(0,\mu_\mathrm{Max})} - 1
= - \frac{p_h}{p} S(0,\mu_\mathrm{Max}) - 1 = -1 < 0.
\]
An analogous computation shows that
\[
S_\mu(0,\mu_\mathrm{Max}) = \frac{1}{p} \int_0^p \frac{U^2}{2}\,\mathrm{d}x > 0,
\]
which completes the proof of the lemma.
\end{proof}


\section{Dynamics near the boundary layer} \label{sec:BoundaryLayer}

In this section, we will prove the existence of a solution to (\ref{ODEPert}) on the interval $[0,r_0]$ for small positive $r_0$.

\begin{lemma} \label{lem:BoundaryLayer}
Assume that Hypotheses~\ref{hyp:Reverser},~\ref{hyp:TrivialEquib}, and~\ref{hyp:Pert} are met. For each compact set $B$ in $\mathrm{Fix}(\mathcal{R})$, there exist $C,R_0,\varepsilon_0>0$ such that for all $\varepsilon\in[0,\varepsilon_0)$, $\mu\in J$, and $u_0\in B$ there exists a unique solution $u=u_\mathrm{bdy}(\cdot;u_0,\mu,\varepsilon)\in C^0([0,2R_0],\mathbb{R}^4)\cap C^1((0,2R_0),\mathbb{R}^4)$ of $(\ref{ODEPert})$ that satisfies the initial condition $u(0)=u_0$. Furthermore, this solution is of the form
\[
u_\mathrm{bdy}(x) = u^0(x) + \bar{u}(x,u_0,\mu,\varepsilon),
\]
where $u^0(x)$ satisfies the unperturbed system (\ref{ODE}) with $u^0(0) = u_0$, and $\bar{u}$ depends smoothly on $(u_0,\mu,\varepsilon)$ with
\[
|\bar{u}(x,u_0,\mu,\varepsilon)|,\,
|\bar{u}_{u_0}(x,u_0,\mu,\varepsilon)|,\,
|\bar{u}_\mu(x,u_0,\mu,\varepsilon)|,\,
|\bar{u}_\varepsilon(x,u_0,\mu,\varepsilon)|
\leq C\varepsilon x
\]
uniformly in $0<x\leq R_0$.
\end{lemma}

\begin{proof}
Let $u^0(x)$ be the solution of $u_x = f(u,\mu)$ with $u^0(0) = u_0 \in \mathrm{Fix}(\mathcal{R})$. Writing $u(x) = u^0(x) + v(x)$, we see that $u(x)$ is a solution to (\ref{ODEPert}) of the form stated in the lemma if and only if $v(x)$ satisfies the nonautonomous initial-value problem
\begin{equation} \label{BndLayer1}
v_x = f(u^0(x)+v,\mu) - f(u^0(x),\mu) + \frac{\varepsilon}{x}g(u^0(x)+v,\mu,\varepsilon), \qquad v(0)=0.
\end{equation}
We write $f(u^0(x)+v,\mu)-f(u^0(x),\mu)=:\tilde{f}(x,v,\mu)v$. Denoting the projection onto $\mathrm{Fix}(-\mathcal{R})$ along $\mathrm{Fix}(\mathcal{R})$ by $\mathcal{P}_\mathcal{R}$, Hypothesis~\ref{hyp:Pert} implies that we have
\[
\frac{1}{x}g(u,\mu,\varepsilon) =: \frac{1}{x}\tilde{g}(u,\mu,\varepsilon)\mathcal{P}_\mathcal{R}u.
\]
Writing $u^0(x)=u_0+x\tilde{u}(x)$ and using that $\mathcal{P}_\mathcal{R}u_0=0$, equation (\ref{BndLayer1}) becomes
\[
v_x = \tilde{f}(x,v,\mu) v + \frac{1}{x}\tilde{g}(u^0(x)+v,\mu,\varepsilon)\mathcal{P}_\mathcal{R} (x \tilde{u}(x)+v).
\]
For $x\in X$, where X is the Banach space defined by
\[
X := \left\{ v\in C^0([0,2R_0],\mathbb{R}^4):\, v(0)=0 \mbox{ and } \|v\|:=\sup_{x\in[0,R_0]}\frac{|v(x)|}{x} <\infty \right\},
\]
we then define a new function $Tv$ by
\begin{eqnarray*}
[Tv](x) & := &
\int_0^x \left[ \tilde{f}(s,v(s),\mu) v(s) + \frac{\varepsilon}{s} \tilde{g}(u^0(s) + v(s),\mu,\varepsilon) \mathcal{P}_\mathcal{R} (s\tilde{u}(s) + v(s)) \right]\,\mathrm{d}s \\ & =: &
\int_0^x \left[ \tilde{f}(s,v(s),\mu) v(s) + \varepsilon \tilde{g}_1(s, v(s),\mu,\varepsilon) + \frac{\varepsilon}{s} \tilde{g}_2(s,v(s),\mu,\varepsilon) v(s) \right]\,\mathrm{d}s.
\end{eqnarray*}
Since fixed points of $T$ are in one-to-one correspondence with solutions of (\ref{BndLayer1}), it suffices to show that, for sufficiently small $R_0,\delta>0$, $T$ maps the ball of radius $\delta$ centered at the origin in $X$ into itself and is a uniform contraction on this ball: these properties are straightforward to verify using the uniform bounds on the smooth functions $\tilde{f}$ and $\tilde{g}_{1,2}$ and their Lipschitz constants in $v$. We omit the details.
\end{proof}



\section{Dynamics near the family of periodic orbits} \label{sec:Fenichel}

The results of the preceding section allow us to restrict the analysis of (\ref{ODEPert}) to the region $x\geq r_0$ for each fixed, but arbitrary, positive value of $r_0$. For each such fixed $r_0>0$, we will construct a local coordinate system akin to Shilnikov variables for the nonautonomous system (\ref{ODEPert}) near the manifold $\mathcal{P}(\mu)$ of periodic orbits that allows us to track solutions as they pass near $\mathcal{P}(\mu)$.

First, note that there exists a closed interval $K_\mathrm{e}$ with $K\subset\mathring{K}_\mathrm{e}$ so that Hypothesis~\ref{hyp:Periodic} holds for all $k\in K_\mathrm{e}$. Our goal is to parametrize the periodic solutions $\gamma(v^c,\mu,v^h)$ by their phase $v^c\in\mathbb{R}$ and the value $v^h$ of the conserved quantity $\mathcal{H}$. We will also use the variables $v^s$ and $v^u$ to parametrize their strong stable and unstable fibers, so that a full neighborhood of the manifold
\[
\mathcal{P}_\mathrm{e}(\mu):=\{\gamma(v^c,\mu,v^h):\, (v^c,v^h)\in\mathbb{R}\times K_\mathrm{e}\}
\]
of periodic orbits is parametrized by $v=(v^c,v^h,v^s,v^u)$. For each $\delta>0$, we define $\mathcal{I}:=[-\delta,\delta]$ and allow $(v^c,v^h,v^s,v^u)$ to vary in the sets
\[
\mathcal{V}_\mathrm{e}:=\mathbb{R}\times K_\mathrm{e}\times\mathcal{I}\times\mathcal{I}, \qquad
\mathcal{V}:=\mathbb{R}\times K\times\mathcal{I}\times\mathcal{I}.
\]
Next, we write (\ref{ODEPert}) as the autonomous system
\begin{equation} \label{e1}
u_x = f(u,\mu) + \varepsilon \rho g(u,\mu,\varepsilon), \qquad
\rho_x = -\rho^2,
\end{equation}
where $\rho=1/x$.

\begin{lemma} \label{lem:FenichelCoords}
Assume that Hypotheses~\ref{hyp:Reverser}, \ref{hyp:Energy}, \ref{hyp:Periodic}, and~\ref{hyp:Pert} are met and recall the constant $R_0$ from Lemma~\ref{lem:BoundaryLayer}, then there are constants $\delta,\varepsilon_1>0$ so that the following is true for each fixed $r_0\in(0,R_0]$. There exist smooth real-valued functions
\[
h^c_{1,2,3},h^e_{1,2,3},h^s,h^u: \mathcal{V}_\mathrm{e}\times J\times[0,\varepsilon_1) \longrightarrow \mathbb{R}, \quad
(v,\mu,\varepsilon) \longmapsto h^c_{1,2,3},h^e_{1,2,3},h^s,h^u(v,\mu,\varepsilon),
\]
the associated differential equation
\begin{equation} \label{Fenichel}
\begin{split}
v^c_x &= 1 + \varepsilon\rho h_1^c(v^c,v^h,\mu,\varepsilon) + \varepsilon^2 \rho^2 h_2^c(\rho,v^c,v^h,\mu,\varepsilon) + h_3^c(\rho,v,\mu,\varepsilon)v^sv^u \\
v^h_x &= \varepsilon \rho h_1^e(v^c,v^h,\mu,\varepsilon) + \varepsilon^2 \rho^2 h_2^e(\rho,v^c,v^h,\mu,\varepsilon) + \varepsilon \rho h^e_3(\rho,v,\mu,\varepsilon)v^sv^u \\
v^s_x &= -[\alpha(\mu,v^h) + h^s(\rho,v,\mu,\varepsilon)]v^s \\
v^u_x &= [\alpha(\mu,v^h) + h^u(\rho,v,\mu,\varepsilon)]v^u \\
\rho_x &= -\rho^2,
\end{split}
\end{equation}
where $x\geq r_0$ and $v\in\mathcal{V}_\mathrm{e}$, and a diffeomorphism from $\mathcal{V}_\mathrm{e}\times J\times[0,\varepsilon_1)$ into a neighborhood of the manifolds $\mathcal{P}_\mathrm{e}(\mu)$ that conjugates (\ref{ODEPert}) and (\ref{Fenichel}) restricted to $\mathcal{V}$ for all $\mu\in J$ and $0\leq\varepsilon<\varepsilon_1$. The functions $h^c_{1,2,3},h^e_{1,2,3},h^s,h^u$ are $p(\mu,v^h)$-periodic in $v^c$, uniformly bounded in $\rho\leq1/r_0$, and globally Lipschitz in $v\in\mathcal{V}_\mathrm{e}$, the functions $h^e_{1,2}$ vanish identically for $v^h\in\partial K_\mathrm{e}$, the functions $h^s$ and $h^u$ vanish identically when $(v^s,v^u,\varepsilon)=0$, and the reverser $\mathcal{R}$ acts on the new coordinates $v$ via $\mathcal{R}(v^c,v^h,v^s,v^u)=(-v^c,v^h,v^u,v^s)$.
\end{lemma}

\begin{proof}
First, we define $v^h:=\mathcal{H}(u,\mu)$ and let $u(x)$ be any solution of $(\ref{e1})$. Using Hypothesis~\ref{hyp:Energy}, we find that
\begin{equation}\label{EnergyFunction}
v^h_x = \nabla_u \mathcal{H}(u,\mu) \cdot u_x
= \nabla_u \mathcal{H}(u,\mu) \cdot \left(f(u,\mu) + \varepsilon \rho g(u,\mu,\varepsilon)\right)
= \varepsilon \rho \nabla_u \mathcal{H}(u,\mu)\cdot g(u,\mu,\varepsilon).
\end{equation}
Next, we set $\varepsilon=0$ and note that the equations for $u$ and $\rho$ in (\ref{e1}) then decouple. We focus first on the equation for $u$ and, using Hypothesis~\ref{hyp:Periodic}(iv), proceed as in \cite{Beck,Fenichel} or \cite[Theorem~4.1.2]{Kuehn} to introduce an invertible coordinate transformation $u=Q(v,\mu)$ defined for $v=(v^c,v^h,v^s,v^u)\in\mathcal{V}_\mathrm{e}$ so that $Q(v^c,v^h,0,0,\mu)=\gamma(v^c,\mu,v^h)$ parametrizes the periodic orbits, the sets $\{v^u=0\}$ and $\{v^s=0\}$ parametrize, for each fixed $(v^c,v^h)$, the strong stable and strong unstable fibers of $\gamma(v^c,\mu,v^h)$, respectively, and the set $\{v^h=h\}$ is equal to $\mathcal{H}^{-1}(h)$. Referring again to \cite{Beck,Fenichel} and extending the coordinate transformation to $(u,\rho)$ via $(u,\rho)=(Q(v,\mu),\rho)$, the vector field in the new coordinates $(v,\rho)$ with $0\leq\rho\leq1/r_0$ is then given by
\begin{equation} \label{FenichelOld}
\begin{split}
v^c_x &= 1 + f^c(v,\mu) v^s v^u +  \varepsilon\rho \tilde{f}^c(\rho,v,\mu,\varepsilon) \\
v^h_x &= \varepsilon\rho \tilde{f}^h(\rho,v,\mu,\varepsilon) \\
v^s_x &= -[\alpha(\mu,v^h) + f^s(v,\mu)] v^s + \varepsilon\rho \tilde{f}^s(\rho,v,\mu,\varepsilon) \\
v^u_x &= [\alpha(\mu,v^h) + f^u(v,\mu)] v^u + \varepsilon\rho \tilde{f}^u(\rho,v,\mu,\varepsilon) \\
\rho_x &= -\rho^2
\end{split}
\end{equation}
for each $\mu\in J$ and $\varepsilon\geq0$, where $f^{s,u}(v,\mu)=0$ when $(v^s,v^u)=0$, and the functions $\tilde{f}^j$ with $j=c,s,u,h$ represent the terms coming from the perturbation $g(Q(v),\mu,\varepsilon)$ with $\tilde{f}^h$ given by (\ref{EnergyFunction}).

When $\varepsilon=0$, the set $\{v^s=v^u=0\}$ is an invariant, normally hyperbolic manifold of (\ref{FenichelOld}). It follows from \cite{Fenichel2,Fenichel} or \cite[Theorem~3.1.4]{Kuehn} that this manifold persists as an invariant, normally hyperbolic manifold of (\ref{FenichelOld}) for all $\varepsilon$ near zero. Since this manifold is given by $(v^s,v^u)=\varepsilon\rho h(v^c,v^h,\rho,\mu,\varepsilon)$ for some smooth function $h$, we can change the $(v^s,v^u)$ coordinates as functions of $(v,\rho)$ so that this three-dimensional manifold is, in the new variables, still given by the set $\{v^s=v^u=0\}$ for all $\varepsilon$ near zero. Similarly, we can straighten out the associated four-dimensional stable and unstable manifolds so that the sets $\{v^u=0\}$, $\{v^s=0\}$, and $\{v^s=v^u=0\}$ are each invariant. The resulting system is of the form (\ref{FenichelOld}) with $(\tilde{f}^s,\tilde{f}^u)=0$ and $(f^s,f^u)$ now depending also on $(\rho,\varepsilon)$.

Since (\ref{FenichelOld}) is autonomous, we know again from \cite{Fenichel2,Fenichel} or \cite[Theorem~3.1.4]{Kuehn} that the four-dimensional invariant stable and unstable manifolds $\{v^u=0\}$ and $\{v^s=0\}$ are each foliated by smooth one-dimensional strong stable and unstable fibers, respectively, that depend smoothly on their base points $(v^c,v^h,0,0,\rho)$. Since these fibers are characterized by exponential decay in forward or backward time towards the solution on the manifold $\{v^s=v^u=0\}$ that passes through the base point, and since we can choose cutoff functions so that the $\rho$ component decays algebraically in both time directions, each fixed fiber must be contained in an appropriate $\{\rho=\mbox{constant}\}$ section. Hence, straightening out these fibers as in \cite[Theorem~4.1.2]{Kuehn}, so that they are given by line segments that do not depend on the base points $(v^c,v^h,0,0,\rho)$, will change only the equations for $(v^c,v^h)$ but not the equation for $\rho$: this change of coordinates then brings (\ref{FenichelOld}) into the normal form (\ref{Fenichel}).

Next, we consider the reverser. The identity $\mathcal{H}(\mathcal{R}u,\mu)=\mathcal{H}(u,\mu)$ implies that $v^h$ remains unchanged under the action of $\mathcal{R}$. It follows from Hypothesis~\ref{hyp:Periodic}(ii) that the action of the reverser on the remaining variables is initially as stated in the lemma, and it is not difficult to check that the action does not change throughout the transformations carried out above. This proves the statements~(i) and~(ii). To establish (iii), we multiply the functions $h^e_{1,2}$ by appropriate cutoff functions so that the products coincide with the original functions for $v^h\in K$ and vanish identically when $v^h\in\partial K_\mathrm{e}$.
\end{proof}

The autonomous vector field (\ref{Fenichel}) for $(v,\rho)$ can be written equivalently as a nonautonomous system for $v$, since we can solve the equation for $\rho$ explicitly to get $\rho(x)=1/x$, which we can then substitute into the remaining equations for $v$. In the remainder of this paper, we will use these two equivalent formulations interchangeably. We define the sets
\[
\mathcal{P}_\mathrm{e}(\mu,\varepsilon) :=\{ v=(v^c,v^h,0,0);\; (v^c,v^h)\in\mathbb{R}\times K_\mathrm{e} \}, \qquad
\mathcal{P}(\mu,\varepsilon) := \{ v=(v^c,v^h,0,0);\; (v^c,v^h)\in\mathbb{R}\times K \}
\]
and note that $\mathcal{P}(\mu,\varepsilon)\subset\mathcal{P}_\mathrm{e}(\mu,\varepsilon)$. Lemma~\ref{lem:FenichelCoords} implies that $\mathcal{P}_\mathrm{e}(\mu,\varepsilon)$ is invariant under the nonautonomous formulation of (\ref{Fenichel}) and that the restriction of this system to $\mathcal{P}_\mathrm{e}(\mu,\varepsilon)$ is given by
\begin{equation} \label{InvariantEqn}
\begin{pmatrix} v^c \\ v^h \end{pmatrix}_x = \begin{pmatrix} 1 \\ 0 \end{pmatrix}
+ \frac{\varepsilon}{x}h_1(v^c,v^h,\mu,\varepsilon) + \frac{\varepsilon^2}{x^2}h_2\left(\frac{1}{x},v^c,v^h,\mu,\varepsilon\right)
=: \begin{pmatrix} 1 \\ 0 \end{pmatrix} + \frac{\varepsilon}{x} F^c(x,v^c,v^h,\mu,\varepsilon).
\end{equation}
We first provide expansions of solutions of (\ref{InvariantEqn}) on $\mathcal{P}_\mathrm{e}(\mu,\varepsilon)$ and show for how long they stay on the smaller manifold $\mathcal{P}(\mu,\varepsilon)$ when they start or end at $v^h=0$.

\begin{lemma}\label{lem:CutInvEqn}
Assume that Hypotheses~\ref{hyp:Reverser}, \ref{hyp:Energy}, \ref{hyp:Periodic}, and~\ref{hyp:Pert} are met. For each fixed $r_0\in(0,1]$ and each sufficiently small $\delta>0$, there are constants $b,C,\varepsilon_1>0$ and a smooth function $g^c(L,\varphi,h,\mu,\varepsilon)$ so that the following is true for each $(\mu,\varepsilon)\in J\times[0,\varepsilon_1)$. 
\begin{compactenum}
\item
For each $(\varphi,h)\in\mathbb{R}\times K_\mathrm{e}$, and $L>r_0$, there exists a unique solution $\Phi(x;r_0,L,\varphi,h,\mu,\varepsilon)$ of (\ref{InvariantEqn}) in $\mathbb{R}\times K_\mathrm{e}$ that satisfies the boundary conditions
\begin{equation}\label{InvariantEqnBCs}
v^c(L) = \varphi, \quad v^h(L)= h
\end{equation}
and lies in $\mathcal{I}\times K_\mathrm{e}$ for $x\in[r_0,L]$. This solution is smooth in $(x,r_0,L,\varphi,h,\mu,\varepsilon)$ and we have
\[
v^c(r_0) = \varphi - L + r_0 + g^c(L,\varphi,h,\mu,\varepsilon)
\]
with $|g^c(L,\varphi,h,\mu,\varepsilon)|+|D_{(\varphi,h_1,\mu)}v^h(r_0)|\leq C\varepsilon\ln L$.
\item
For each $(\varphi,h)\in\mathbb{R}\times U_\delta(0)$ and $r_0<L\leq\mathrm{e}^{b/\varepsilon}$, the solution $(v^c(x),v^h(x))$ of (\ref{InvariantEqn})--(\ref{InvariantEqnBCs}) satisfies $v^h(x)\in K$ for all $x\in[r_0,L]$.
\end{compactenum}
\end{lemma}

\begin{proof}
To prove (i), we note that existence, uniqueness, and smoothness of the solution follows since $\mathcal{I}\times K_\mathrm{e}$ is invariant under (\ref{InvariantEqn}). It therefore remains to estimate $v^c(r_0)$. We write $v^c(x)=\varphi-L+x+\tilde{v}^c(x)$ so that the boundary condition $v^c(L)=\varphi$ becomes $\tilde{v}^c(L)=0$, and $\tilde{v}^c(r_0)=g^c(L,\varepsilon,h,\mu,\varepsilon)$ is then given implicitly by
\[
\tilde{v}^c(r_0) = \varepsilon \int_L^{r_0}\frac{h_1^c(v^c(s),v^h(s),\mu,\varepsilon)}{s}\,\mathrm{d}s + \varepsilon^2\int_L^{r_0}\frac{h_2^c(1/s,v^c(s),v^h(s),\mu,\varepsilon)}{s^2}\,\mathrm{d}s.
\]
Bounding $|h^c_{1,2}|$ by a uniform constant $C_0>0$, we obtain
\[
|\tilde{v}^c(r_0)| = |g^c(L,\mu,\varepsilon)| \leq \varepsilon C_0 (|\ln L|+|\ln r_0|) \leq \varepsilon C \ln L.
\]
The bounds on $|v^h(r_0)|$ and its derivatives are handled in an identical manner. For (ii), it suffices to find conditions that guarantee that $v^h(x)\in K_\mathrm{e}$ for all $x\in[r_0,L]$ whenever $v^h(L)=h\in U_\delta(0)$. Since $v^h(x)$ satisfies
\[
v^h(x) = h + \varepsilon \int_L^x\frac{h_1^e(v^c(s),v^h(s),\mu,\varepsilon)}{s}\,\mathrm{d}s
+ \varepsilon^2\int_L^x\frac{h_2^e(1/s,v^c(s),v^h(s),\mu,\varepsilon)}{s^2}\,\mathrm{d}s,
\]
and we can bound $|h^e_{1,2}|$ by a uniform constant $C_0$, we find that
\[
|v^h(x)| \leq |h| + \varepsilon C_0 (|\ln L|+|\ln r_0|)
\]
for $r_0\leq x\leq L$. If $K=[-k,k]$, then setting $\varepsilon_1=\frac{k-\delta}{2C_0|\ln r_0|}$ and $b:=\frac{k-\delta}{2C_0}$ guarantees that $|v^h(x)|\leq k$ for $r_0\leq x\leq L$ for $L\leq\mathrm{e}^{b/\varepsilon}$ as claimed.
\end{proof}

Next, we use the full equation (\ref{Fenichel}) to track solutions as they evolve in $\mathcal{V}_\mathrm{e}$ near $\mathcal{P}_\mathrm{e}(\mu,\varepsilon)$.

\begin{proposition} \label{prop:Fenichel}
Assume that Hypotheses~\ref{hyp:Reverser}, \ref{hyp:Energy}, \ref{hyp:Periodic}, and~\ref{hyp:Pert} are met. There exist constants $\eta,L_0,M>0$ such that, for each fixed $0<r_0\leq1$, there exists an $\varepsilon_2>0$ such that the following holds: pick $0\leq\varepsilon\leq\varepsilon_2$ and $L\geq L_0$, and let $\Phi(x;r_0,L,\varphi,h,\mu,\varepsilon)$ be as in Lemma~\ref{lem:CutInvEqn}, then, for each $a^s\in\mathcal{I}$, there exists a unique solution $v(x)=v(x;r_0,L,a^s,\varphi,h,\mu,\varepsilon)\in\mathcal{V}_\mathrm{e}$ to (\ref{Fenichel}) defined for $x \in [r_0,L]$ so that
\[
v^c(L) = \varphi,\qquad
v^h(L) = h, \qquad
v^s(r_0) = a^s,\qquad
v^u(L) =\delta.
\]
Furthermore, this solution satisfies
\begin{equation} \label{FenichelEstimates}
|(v^c(x),v^h(x))^T - \Phi(x;r_0,L,\varphi,h,\mu,\varepsilon)| \leq M \mathrm{e}^{-\eta L}, \quad
|v^s(x)| \leq M \mathrm{e}^{-\eta x}, \quad
|v^u(x)| \leq M \mathrm{e}^{\eta (x-L)}
\end{equation}
for all $x\in[r_0,L]$, the solution $v(x)$ is smooth in $(x,r_0,L,a^s,\varphi,h,\mu,\varepsilon)$, and the bounds (\ref{FenichelEstimates}) also hold for these derivatives.
\end{proposition}

\begin{proof}
We will show that the assumptions of \cite[Theorem~2.2]{Schecter} are satisfied: our statements then follow directly from this theorem. Note that restricting to $x\geq r_0$ and choosing $0\leq\varepsilon_2\leq r_0^2$ guarantees that the right-hand side of (\ref{Fenichel}) is bounded uniformly in $x$. It remains to establish appropriate exponential bounds for solutions of the linearized dynamics of (\ref{Fenichel}). Linearizing (\ref{Fenichel}) along the solution $(v^c,v^h,v^s,v^u)=(\Phi^c(x),\Phi^h(x),0,0)$, where $(\Phi^c(x),\Phi^h(x))^T = \Phi(x;r_0,L,\mu,\varepsilon)$ satisfies (\ref{InvariantEqn})--(\ref{InvariantEqnBCs}) on $[r_0,L]$, we arrive at the linear system
\begin{subequations}
\begin{align}
w_x &= \frac{\varepsilon}{x} D_{(v^c,v^h)}F^c(x,\Phi^c(x),\Phi^h(x),\varepsilon)c, \label{LinearizedCentre} \\
v^s_x &= -[\alpha(\Phi^h(x))+ h^s(x,\Phi^c(x),0,0,\Phi^h(x),\varepsilon)]v^s \label{LinearizedStable} \\
v^u_x &= [\alpha(\Phi^h(x))+ h^u(x,\Phi^c(x),0,0,\Phi^h(x),\varepsilon)]v^u \label{LinearizedUnstable}
\end{align}
\end{subequations}
where $w = (v^c,v^h)^T$; note that we have suppressed the dependence on $\mu\in J$ for notational convenience. Lemma~\ref{lem:FenichelCoords} implies that $(h^s,h^u)$ vanish uniformly when $\varepsilon = 0$, and Hypothesis~\ref{hyp:Periodic} implies that $\alpha(\mu,v^h)$ is bounded away from zero uniformly in $(\mu,v^h) \in J \times K$. Hence, for $\varepsilon > 0$ taken sufficiently small, the right-hand sides of (\ref{LinearizedStable}) and (\ref{LinearizedUnstable}) are uniformly bounded away from zero, and we conclude that there are constants $\eta^s, \eta^u > 0$ and $M_0 > 0$ such that the solution operators $\Psi^s(x,s),\Psi^u(x,s)$ of (\ref{LinearizedStable}) and (\ref{LinearizedUnstable}), respectively, satisfy
\[
|\Psi^s(x,y)| \leq M_0 \mathrm{e}^{-\eta^s (x-y)} \mbox{ and } |\Psi^u(y,x)| \leq M_0 \mathrm{e}^{\eta^u (y-x)}
\]
for $r_0\leq y\leq x\leq L$ and $\varepsilon>0$ sufficiently small. Furthermore, we may take $\eta^u>0$ so that $-\eta^s+\eta^u < 0$.

We now turn to (\ref{LinearizedCentre}). Lemma~\ref{lem:FenichelCoords} shows that there is a constant $C>0$ with $|D_{(v^c,v^h)}F^c(x,\Phi^c(x),\Phi^h(x),\varepsilon)|\leq C$ uniformly in all arguments, and we conclude that
\[
|w_x| \leq \frac{\varepsilon C}{x} |w| \leq \varepsilon^{\frac12} C |w|
\]
provided we pick $0\leq\varepsilon_2\leq r_0^2$. Denoting the solution operator to (\ref{LinearizedCentre}) by $\Psi^c(x,y)$, we have
\[
|\Psi^c(x,y)| \leq \mathrm{e}^{\sqrt{\varepsilon}C(x-y)}
\]
for all $x,y\in [r_0,L]$, independently of $r_0$ and $L$, which verifies \cite[Hypothesis~(E1) in Theorem~2.2]{Schecter}. Finally, taking $\varepsilon > 0$ and sufficiently small, we can guarantee that
\[
\eta^s + \eta^u + \sqrt{\varepsilon}K_2 < 0 < \eta^u - \sqrt{\varepsilon}C,
\]
which verifies \cite[Hypothesis~(D2) in Theorem~2.2]{Schecter}. We have now verified all hypotheses required to apply \cite[Theorem~2.2]{Schecter}, which proves the statement of Proposition~\ref{prop:Fenichel}.
\end{proof}


\section{Dynamics near the stable manifold of the homogeneous state} \label{sec:StableManifold}

Hypothesis~\ref{hyp:Pert} implies that $u(x)=0$ satisfies (\ref{ODEPert}) for all $\varepsilon\geq0$, and we now describe the set of solutions of (\ref{ODEPert}) that converge to zero as $x\to\infty$ for $\varepsilon\geq0$. For each $L\geq1$, we define
\[
W^s_L(0,\mu,\varepsilon) := \{u_0 \in \mathbb{R}^4:\, u(x) \mbox{ satisfies } (\ref{ODEPert}) \mbox{ with } u(L) = u_0, \mbox{ and } u(x)\to0 \mbox{ as } x\to\infty\}
\]
as the section of the stable manifold of the trivial solution at $x=L$. First, we show that $W^s_L(0,\mu,\varepsilon)$ is a regular perturbation of $W^s(0,\mu)$ in $\varepsilon$.

\begin{lemma} \label{lem:TrivialStable}
Assume that Hypotheses~\ref{hyp:TrivialEquib} and~\ref{hyp:Pert} are met, then for each $L\geq1$ the set $W^s_L(0,\mu,\varepsilon)$ is a two-dimensional manifold that is locally $\mathcal{O}(\varepsilon/L)$ $C^1$-close to $W^s(0,\mu)$ uniformly in $\mu$.
\end{lemma}

Lemma~\ref{lem:TrivialStable} follows directly from the uniform contraction mapping principle, and we omit the details. Next, we use this lemma to provide a parametrization of the stable manifold in the Shilnikov variables that were introduced in the preceding section.

\begin{lemma} \label{lem:GFunction}
Assume that Hypotheses~\ref{hyp:Reverser}-\ref{hyp:Pert} are met and define
\[
\Sigma_\mathrm{out} = \mathbb{R}\times K_\mathrm{e}\times\mathcal{I}\times\{v^u =\delta\}.
\]
There exist $\varepsilon_0>0$ and smooth real-valued functions $z^\Gamma,z^h$ so that the following is true for each $\varepsilon\in[0,\varepsilon_0]$ and $L\geq1$: A point $(v^c,v^h,v^s,\delta)$ with $|v^s|,|v^h|<\varepsilon_0$ lies in $W^s_L(0,\mu,\varepsilon)\cap\Sigma_\mathrm{out}$ if and only if there exists $(\varphi^*,\mu^*)\in\Gamma$ such that
\begin{eqnarray*}
v^c & = & \varphi^* + z^\Gamma(L,\varphi^*,v^s,\mu^*,\varepsilon) G_\varphi(\varphi^*,0,\mu^*) \\
\mu & = & \mu^* + z^\Gamma(L,\varphi^*,v^s,\mu^*,\varepsilon) G_\mu(\varphi^*,0,\mu^*) \\
v^h & = & z^h(L,\varphi^*,v^s,\mu^*,\varepsilon),
\end{eqnarray*}
where the function $G$ was defined in Hypothesis~\ref{hyp:StableIntersection}. Furthermore, the functions $z^\Gamma$ and $z^h$ are bounded uniformly in $L\geq1$, independent of $L$ when $\varepsilon=0$, and satisfy
\[
z^\Gamma(\cdot,\varphi^*,0,\mu^*,0) = 0, \qquad
z^h(\cdot,\varphi^*,v^s,\mu^*,0) = 0.
\]
\end{lemma}

\begin{proof}
Define the subset of $\Sigma_\mathrm{out}\times J\cong\mathbb{R}\times K_\mathrm{e}\times\mathcal{I}\times J$ given by
\[
\tilde{W}^s_L(0,\varepsilon) := \bigcup\limits_{\mu \in J} \left(W^s_L(0,\mu,\varepsilon)\cap \{v^u = \delta\}\right) \times \{\mu\}.
\]
Note that $\tilde{W}^s_L(0,0)$ is independent of $L$ since (\ref{ODEPert}) is autonomous when $\varepsilon = 0$. Hypothesis~\ref{hyp:StableIntersection} implies that
\[
\Gamma = \tilde{W}^s_L(0,0) \cap (\mathbb{R} \times \{0\} \times \{0\} \times J),
\]
that the function $\mathcal{G}:\mathbb{R}\times K_\mathrm{e}\times\mathcal{I}\times J \to \mathbb{R}^2$ given by
\[
\mathcal{G}(v^c,v^h,v^s,\mu) = (G(v^c,v^s,\mu),v^h)
\]
satisfies $\mathcal{G}^{-1}(0) = \tilde{W}^s_L(0,0)$, and that $D_{(v^c,v^h,\mu)} \mathcal{G}(v^c,0,0,\mu)\in\mathbb{R}^{2\times 3}$ has full rank for all $(v^c,\mu)\in\Gamma$. Therefore, we may use the columns of $D_{(v^c,v^h,\mu)} \mathcal{G}(v^c,0,0,\mu)$ to define normal vectors to $\Gamma$ as a subset of $\tilde{W}^s_L(0,\varepsilon)$ in the space $\Sigma_\mathrm{out} \times J$ to construct the functions $z^\Gamma$ and $z^h$ that represent the displacement from $\Gamma$ along these normal vectors with respect to small perturbations in $v^s$ and $\varepsilon$. Smoothness and boundedness with respect to $L\geq1$ follow from Lemma~\ref{lem:TrivialStable}. The property that $z^h(\varphi^*,v^s,\mu^*,0) = 0$ follows from the fact that $W^s(0,\mu)\subset\mathcal{H}^{-1}(0)$ for all $\mu \in J$. This completes the proof.
\end{proof}


\section{Construction of radial pulses} \label{sec:Matching}

To construct radial pulses, we need to match the solution segments we obtained in the preceding sections on the spatial intervals $[0,r_0]$, $[r_0,L]$, and $[L,\infty)$. Recall that $u(x)$ is a radial pulse of (\ref{ODEPert}) if
\[
u(0)\in\mathrm{Fix}(\mathcal{R}), \qquad
u(x)\in\mathcal{V} \mbox{ for } 0\leq x\leq L, \qquad
u(L)\in W_L^s(0,\mu,\varepsilon)\cap\Sigma_\mathrm{out}.
\]
In this section, we will construct radial pseudo-pulses which, by definition, satisfy
\[
u(0) \in \mathrm{Fix}(\mathcal{R}), \qquad
u(x) \in \mathcal{V}_\mathrm{e} \mbox{ for } 0\leq x\leq L, \qquad
u(L) \in W_L^s(0,\mu,\varepsilon) \cap \Sigma_\mathrm{out}.
\]
Note that radial pseudo-pulses correspond to radial pulses of (\ref{ODEPert}) only when the values of $v^h$ for $0\leq x\leq L$ lie in $K$ rather than in the larger set $K_\mathrm{e}$. Before stating our result, we recall that $\Gamma_\mathrm{lift}\subset\mathbb{R}\times\mathring{J}$ is the preimage of the set $\Gamma$ defined in (\ref{Gamma}) under the natural covering map from $\mathbb{R}\times\mathring{J}$ to $S^1\times\mathring{J}$.

\begin{theorem} \label{thm:epsilonsmall2}
Assume that Hypotheses~\ref{hyp:Reverser}-\ref{hyp:Pert} are met. There exist constants $\varepsilon_0,\eta,L_*>0$ and sets $\Gamma_\mathrm{pulse}^{\varphi_0,\varepsilon} \subset (L_*,\infty) \times J$ defined for $\varphi_0\in\{0,\pi\}$ and $\varepsilon \in [0,\varepsilon_0]$ so that the following is true:
\begin{compactenum}
\item A radial pseudo-pulse with plateau length $L$ exists if and only if $(L,\mu)\in\Gamma_\mathrm{pulse}^{\varphi_0,\varepsilon}$ for $\varphi_0 = 0$ or $\varphi_0 = \pi$.
\item There exists a smooth function $\tilde{g}_\mathrm{lift}(L,\mu,\varepsilon) = \mathcal{O}(\varepsilon\ln L)$ such that for each fixed $\varphi_0\in\{0,\pi\}$ and $\varepsilon \in [0,\varepsilon_0]$ the one-dimensional manifolds
\[
\tilde{\Gamma}_\mathrm{lift}^{\varphi_0,\varepsilon} := \{(L-\tilde{g}_\mathrm{lift}(L,\mu,\varepsilon)+\varphi_0,\mu) :\ (L,\mu) \in \Gamma_\mathrm{lift}\cap((L_*\infty) \times J)\}
\]
and $\Gamma_\mathrm{pulse}^{\varphi_0,\varepsilon}$ are $\mathcal{O}(\mathrm{e}^{-\eta L})$-close to each other in the $C^0$-sense near any point $(L,\mu) \in \tilde{\Gamma}_\mathrm{lift}^{\varphi_0,\varepsilon}$.
\end{compactenum}
\end{theorem}

Before we give the proof of this theorem, we show how Theorem~\ref{thm:epsilonsmall} follows from it.

\begin{proof}[Proof of Theorem~\ref{thm:epsilonsmall}]
We need to verify that the radial pseudo-pulses constructed in Theorem~\ref{thm:epsilonsmall2} are radial pulses (meaning that their $v^h$ components stay in $K$ for all $0\leq x\leq L$) under the restrictions on $L$ assumed in Theorem~\ref{thm:epsilonsmall}. This follows from the estimates given in Lemma~\ref{lem:CutInvEqn}(ii) and Proposition~\ref{prop:Fenichel} on the evolution of (\ref{Fenichel}) on $\mathcal{P}_\mathrm{e}(\mu,\varepsilon)$.
\end{proof}

Next, we give the proof of Theorem~\ref{thm:epsilonsmall2}.

\begin{proof}[Proof of Theorem~\ref{thm:epsilonsmall2}]
We need to match (i) the boundary-layer solution $u_\mathrm{b}(x)$ from Lemma~\ref{lem:BoundaryLayer} and the plateau solution from Proposition~\ref{prop:Fenichel} at $x=r_0$ and (ii) the plateau solution and the stable manifold $W^s_L(0,\mu,\varepsilon)$ from Lemma~\ref{lem:GFunction} at $x=L$ in $\Sigma_\mathrm{out}$. We will accomplish this in several steps.

\paragraph{Boundary-layer solution.}
We first transform the boundary-layer solution $u_\mathrm{b}(x)$ obtained in Lemma~\ref{lem:BoundaryLayer} into the Fenichel coordinates derived in Lemma~\ref{lem:FenichelCoords}. Doing so, we see that for each $0<r_0\ll1$, $u_0\in\mathrm{Fix}(\mathcal{R})$, $\varphi_0\in\{0,\pi\}$, $a_0\in\mathcal{I}$, and $h_0 \in K$ we can set $u_0=(\varphi_0,h_0,a_0,a_0)\in\mathrm{Fix}(\mathcal{R})$ and write the corresponding boundary-layer solution $u_\mathrm{b}(x)$ as
\[
v_\mathrm{b}(r_0,u_0,\mu,\varepsilon) = (\varphi_0,h_0,a_0,a_0) + (v_\mathrm{b}^c,v_\mathrm{b}^h,v_\mathrm{b}^s,v_\mathrm{b}^u)(r_0,u_0,\mu,\varepsilon),
\]
where $v_\mathrm{b}^j(0,u_0,\mu,\varepsilon) = 0$ for all $j = c,h,s,u$ and $v_\mathrm{b}(0,u_0,\mu,\varepsilon) = u_0 \in \mathrm{Fix}(\mathcal{R})$. Setting $\varepsilon=0$, we obtain
\[
v_\mathrm{b}(r_0,u_0,\mu,0) = (\varphi_0+r_0,h_0,a_0,a_0) + \mathcal{O}(r_0a_0),
\]
where the $\mathcal{O}(r_0a_0)$ term is independent of $h_0$: indeed, $v_\mathrm{b}$ depends smoothly on $x$, and $a_0=0$ restricts the solution to the cylinder of periodic orbits whose $v^h$-component does not change. Smoothness in $\varepsilon$ then shows that
\begin{equation}\label{v_bExpansion}
v_\mathrm{b}(r_0,u_0,\mu,\varepsilon) = (\varphi_0+r_0,h_0,a_0,a_0) + \mathcal{O}(r_0a_0) + \mathcal{O}(\varepsilon).
\end{equation}

\paragraph{Plateau solution.}
Lemma~\ref{lem:CutInvEqn}(i) shows that for each $\varphi_1\in\mathbb{R}$ and $h_1\in K$ there is a unique solution
\[
(\Phi^c,\Phi^h)(x) = (\Phi^c,\Phi^h)(x;r_0,L,\varphi_1,h_1,\mu,\varepsilon)
\]
with
\[
\Phi^c(L;r_0,L,\varphi_1,h_1,\mu,\varepsilon) = \varphi_1, \qquad
\Phi^h(L;r_0,L,\varphi_1,h_1,\mu,\varepsilon) = h_1
\]
on the invariant manifold $\mathcal{P}_\mathrm{e}(\mu,\varepsilon)$. For each $b_0$, Proposition~\ref{prop:Fenichel} then yields a plateau solution $v_\mathrm{f}(x;r_0,L,b_0,\varphi_1,h_1,\mu,\varepsilon)$, which is smooth in its arguments and satisfies
\begin{equation}\label{FenMatching}
\begin{split}
v_\mathrm{f}(r_0;r_0,L,b_0,\varphi_1,h_1,\mu,\varepsilon) &=
(\Phi^c(r_0)+\mathcal{O}(\mathrm{e}^{-\eta L}), \Phi^h(r_0)+\mathcal{O}(\mathrm{e}^{-\eta L}), b_0, \mathcal{O}(\mathrm{e}^{-\eta L})), \\
v_\mathrm{f}(L;r_0,L,b_0,\varphi_1,h_1,\mu,\varepsilon) &= (\varphi_1, h_1, \mathcal{O}(\mathrm{e}^{-\eta L}),\delta).
\end{split}
\end{equation}
Note that the error estimates in (\ref{FenMatching}) can be differentiated and hold for all derivatives with respect to the arguments of $v_\mathrm{f}$.

\paragraph{Matching solutions at $\mathbf{x=r_0}$.}
Using the solutions introduced above, condition~(i) becomes $v_\mathrm{b}(r_0)=v_\mathrm{f}(r_0)$, which we write in coordinates as
\begin{subequations}
\begin{align}
\varphi_0 + r_0 + \mathcal{O}(r_0a_0 + \varepsilon) &=
\Phi^c(r_0;r_0,L,\varphi_1,h_1,\mu,\varepsilon) + \mathcal{O}(\mathrm{e}^{-\eta L}) \label{r_0Matching1} \\
a_0 + \mathcal{O}(r_0a_0 + \varepsilon) &= b_0 + \mathcal{O}(\mathrm{e}^{-\eta L}) \label{r_0Matching2} \\
a_0 + \mathcal{O}(r_0a_0 + \varepsilon) &= \mathcal{O}(\mathrm{e}^{-\eta L}) \label{r_0Matching3} \\
h_0 + \mathcal{O}(r_0a_0 + \varepsilon) &=
\Phi^h(r_0;r_0,L,\varphi_1,h_1,\mu,\varepsilon) + \mathcal{O}(\mathrm{e}^{-\eta L}) \label{r_0Matching4}
\end{align}
\end{subequations}
using the expansion for $v_\mathrm{b}(r_0)$ stated in (\ref{v_bExpansion}). We focus first on (\ref{r_0Matching2})-(\ref{r_0Matching4}) and will return to (\ref{r_0Matching1}) at the end of the proof. Define a smooth function $F_0:\mathbb{R}^3 \to \mathbb{R}^3$ that depends on the variables $(a_0,b_0,h_0)$ and parameters $(r_0,L,\varphi_0,\varphi_1,h_1,\mu,\varepsilon)$ so that its roots correspond to solutions of (\ref{r_0Matching2})-(\ref{r_0Matching4}). Note that $F_0$ has the expansion
\[
F_0(a_0,b_0,h_0) = \begin{pmatrix}
a_0 - b_0 \\ a_0 \\ h_0 - \Phi^h(r_0;r_0,L,\varphi_1,h_1,\mu,\varepsilon)
\end{pmatrix} + \mathcal{O}(r_0a_0 + \mathrm{e}^{-\eta L} + \varepsilon).
\]
To find its roots, we will use the following result that we state without proof.

\begin{lemma}\label{thm:Roots}
If $F:\mathbb{R}^n \to \mathbb{R}^n$ is smooth and there are constants $0<\kappa<1$ and $\rho>0$, a vector $w_0\in\mathbb{R}^n$, and an invertible matrix $A\in\mathbb{R}^{n\times n}$ so that
\begin{compactenum}
\item $\|1 - A^{-1}DF(w)\| \leq \kappa$ for all $w \in B_\rho(w_0)$, and
\item $\|A^{-1}F(w_0)\|\leq (1-\kappa)\rho$,
\end{compactenum}
then $F$ has a unique root $w_*$ in $B_\rho(w_0)$, and this root satisfies $|w_* - w_0| \leq \frac{1}{1 - \kappa}\|A^{-1}F(w_0)\|$.
\end{lemma}

Using the notation of Lemma~\ref{thm:Roots}, let $w_0= (0,0,\Phi^h(r_0;r_0,L,b_0,\varphi_1,h_1,\mu,\varepsilon))$ and define $A_0$ to be the invertible matrix
\[
A_0 = \begin{bmatrix} 1 & -1 & 0 \\ 0 & 1 & 0 \\ 0 & 0 & 1 \end{bmatrix}.
\]
We then have $F_0(w_0)=\mathcal{O}(\mathrm{e}^{-\eta L}+\varepsilon)$ and therefore $\|A_0^{-1}F_0(w_0)\|=\mathcal{O}(\mathrm{e}^{-\eta L}+\varepsilon)$. Note furthermore that
\[
\|1 - A_0^{-1}DF_0(w)\| = \mathcal{O}(r_0 + \mathrm{e}^{-\eta L} + \varepsilon)
\]
for all $w$ since $\Phi^h(r_0;r_0,L,\varphi_1,h_1,\mu,\varepsilon)$ is independent of $(a_0,b_0,h_0)$. Hence, choosing $\kappa=\frac{1}{2}$, $\rho=1$, and $0<r_0\ll1$ sufficiently small, Lemma~\ref{thm:Roots} shows that there is a unique function $(a_0^*,b_0^*,h_0^*)(L,\varphi_0,\varphi_1,h_1,\mu,\varepsilon)$ that is defined for all $0<\varepsilon\ll1$, $L\gg1$, $\varphi_0\in\{0,\pi\}$, and arbitrary $(\varphi_1,h_1,\mu)$, corresponds to roots of $F_0$, is smooth in the arguments $(L,\varphi_1,h_1,\mu,\varepsilon)$, and has the expansion
\[
(a_0^*,b_0^*,h_0^*)(L,\varphi_0,\varphi_1,h_1,\mu,\varepsilon) = (0,0,\Phi^h(r_0;r_0,L,\varphi_1,h_1,\mu,\varepsilon)) + \mathcal{O}(\mathrm{e}^{-\eta L} + \varepsilon).
\]
Moreover, recalling from Lemma~\ref{lem:CutInvEqn}(i) that
\[
D_{(\varphi_1,\mu,h_1)} \Phi^h(r_0;r_0,L,\varphi_1,h_1,\mu,\varepsilon) = \mathcal{O}(\varepsilon\ln L),
\]
we find that
\begin{equation} \label{epslnDerivative}
D_{(\varphi_1,\mu,h_1)}(a_0^*,b_0^*,h_0^*)(L,\varphi_0,\varphi_1,h_1,\mu,\varepsilon) = \mathcal{O}(\mathrm{e}^{-\eta L} + \varepsilon + \varepsilon\ln L).
\end{equation}

\paragraph{Matching solutions at $\mathbf{x=L}$.}
Next, we consider condition~(ii), which requires that $v_\mathrm{f}(L)\in W^s_L(0,\mu,\varepsilon)\cap\Sigma_\mathrm{out}$. Since we will rely on the characterization of $W^s_L(0,\mu,\varepsilon)$ given in Lemma~\ref{lem:GFunction}, we need to parametrize the set $\Gamma$: we shall follow the construction introduced in \cite{Aougab}. If $\Gamma$ consists of 0-loops, we parametrize each loop by $2\pi$-periodic functions $(\tilde{\varphi}(s),\tilde{\mu}(s))$ with $0\leq s\leq2\pi$ so that
\[
\Gamma_\mathrm{lift} = \{(\tilde{\varphi}(s) + 2\pi j,\tilde{\mu}(s))\ :\ 0\leq s\leq 2\pi, j\in \mathbb{N}\}.
\]
If $\Gamma$ is a 1-loop, we parametrize $\Gamma_\mathrm{lift}$ by a curve $(\tilde{\varphi}(s),\tilde{\mu}(s))$ with $s\geq0$, where $(\tilde{\varphi}(s+2\pi), \tilde{\mu}(s+2\pi))=(\tilde{\varphi}(s)+2\pi,\tilde{\mu}(s))$.

Using Lemma~\ref{lem:GFunction} and the fact that $(v_\mathrm{f}^c,v_\mathrm{f}^h)(L;r_0,L,\varphi_1,h_1,\mu,\varepsilon)=(\varphi_1,h_1)$, we see that the condition $v_\mathrm{f}(L)\in W^s_L(0,\mu,\varepsilon)\cap\Sigma_\mathrm{out}$ is equivalent to the system
\begin{equation} \label{LMatching1}
\begin{split}
\varphi_1 &= \tilde{\varphi}(s) +
z^\Gamma(L,\tilde{\varphi}(s),v^s_\mathrm{f}(L;r_0,L,b_0^*,\varphi_1,h_1,\mu,\varepsilon),\tilde{\mu}(s),\varepsilon)
G_\varphi(\tilde{\varphi}(s),0,\tilde{\mu}(s)), \\
\mu &= \tilde{\mu}(s) +
z^\Gamma(L,\tilde{\varphi}(s),v^s_\mathrm{f}(L;r_0,L,b_0^*,\varphi_1,h_1,\mu,\varepsilon),\tilde{\mu}(s),\varepsilon)
G_\mu(\tilde{\varphi}(s),0,\tilde{\mu}(s)), \\
h_1 &= z^h(L,\tilde{\varphi}(s),v^s_\mathrm{f}(L;r_0,L,b_0^*,\varphi_1,h_1,\mu,\varepsilon),\tilde{\mu}(s),\varepsilon),
\end{split}
\end{equation}
where $b_0^*=b_0^*(L,\varphi_0,\varphi_1,h_1,\mu,\varepsilon)$. Using that $(z^\Gamma,z^h)$ vanish when $(v^s_\mathrm{f},\varepsilon)=0$ and that
\[
|v^s_\mathrm{f}(L;r_0,L,b_0^*,\varphi_1,h_1,\mu,\varepsilon)| \leq M\mathrm{e}^{-\eta L}
\]
by Proposition~\ref{prop:Fenichel}, we conclude that (\ref{LMatching1}) can be written
\begin{equation} \label{LMatching2}
\begin{split}
\varphi_1 - \tilde{\varphi}(s) + \mathcal{O}(\mathrm{e}^{-\eta L} + \varepsilon) & = 0 \\
\mu - \tilde{\mu}(s) + \mathcal{O}(\mathrm{e}^{-\eta L} +\varepsilon) & = 0 \\
h_1 + \mathcal{O}(\mathrm{e}^{-\eta L} + \varepsilon) & = 0.
\end{split}
\end{equation}
Using Proposition~\ref{prop:Fenichel} and the estimate (\ref{epslnDerivative}), we see that the estimates for the remainder terms in (\ref{LMatching2}) hold also for their derivatives with respect to $(\varphi_1,\mu,h_1)$. Thus, (\ref{LMatching2}) is of the form
\[
F_1(\varphi_1,\mu,h_1) = \begin{pmatrix} \varphi_1 - \tilde{\varphi}(s) \\ \mu - \tilde{\mu}(s) \\ h_1 \end{pmatrix} + \mathcal{O}(\mathrm{e}^{-\eta L} + \varepsilon) = 0,
\]
where we consider $(L,\varphi_0,\varepsilon)$ as parameters. As before, we will use Lemma~\ref{thm:Roots} to find solutions to (\ref{LMatching2}) by characterizing roots of $F_1$. We take $w_1=(\tilde{\varphi}(s),\tilde{\mu}(s),0)$ and let $A_1$ be the $3\times3$ identity matrix, which immediately gives $\|A_1^{-1}F_1(w_1)\|=\mathcal{O}(\mathrm{e}^{-\eta L}+\varepsilon)$ and $\|1-A_1^{-1}DF_1(w)\|=\mathcal{O}(\mathrm{e}^{-\eta L}+\varepsilon)$. Applying Lemma~\ref{thm:Roots} with $\kappa=\frac{1}{2}$ and $\rho=1$ shows that roots of $F_1$ are given by a unique function $(\varphi_1^*,\mu^*,h_1^*)(L,\varphi_0,\varepsilon,s)$, which is defined for all $s$, $0<\varepsilon\ll1$, $L\gg1$, and $\varphi_0\in\{0,\pi\}$, depends smoothly on its arguments $(L,\varepsilon,s)$, and has the expansion
\[
(\varphi_1^*,\mu^*,h_1^*)(L,\varphi_0,\varepsilon,s) =
(\tilde{\varphi}(s),\tilde{\mu}(s),0) + \mathcal{O}(\mathrm{e}^{-\eta L}+\varepsilon).
\]

In summary, upon evaluating $(a_0^*,b_0^*,h_0^*)$ at $(\varphi_1^*,\mu^*,h_1^*)(L,\varphi_0,\varepsilon,s)$, we see that the solutions to (\ref{r_0Matching2})--(\ref{r_0Matching4}) and (\ref{LMatching1}) are of the form
\begin{equation} \label{MatchingSoln}
\begin{split}
a_0^*(L,\varphi_0,\varepsilon,s) &= \mathcal{O}(\mathrm{e}^{-\eta L} + \varepsilon), \\
b_0^*(L,\varphi_0,\varepsilon,s) &= \mathcal{O}(\mathrm{e}^{-\eta L} + \varepsilon), \\
h_0^*(L,\varphi_0,\varepsilon,s) &= \Phi^h(r_0;r_0,L,\varphi_1(s),0,\mu(s),\varepsilon) + \mathcal{O}(\mathrm{e}^{-\eta L} + \varepsilon + \varepsilon^2\ln L), \\
\varphi_1^*(L,\varphi_0,\varepsilon,s) &= \tilde{\varphi}(s) + \mathcal{O}(\mathrm{e}^{-\eta L} + \varepsilon), \\
\mu^*(L,\varphi_0,\varepsilon,s) &= \tilde{\mu}(s) + \mathcal{O}(\mathrm{e}^{-\eta L} + \varepsilon), \\
h_1^*(L,\varphi_0,\varepsilon,s) &= \mathcal{O}(\mathrm{e}^{-\eta L} + \varepsilon).
\end{split}
\end{equation}
In particular, the only remaining free variables are $(L,\varphi_0,\varepsilon,s)$. We now return to the remaining equation (\ref{r_0Matching1}).

\paragraph{Matching the phase at $\mathbf{x=r_0}$.}
It remains to solve (\ref{r_0Matching1}), which, upon substituting the expressions  (\ref{MatchingSoln}), becomes
\begin{equation}\label{LSoln}
\varphi_0 + r_0 + \mathcal{O}(r_0a_0^*(L,\varphi_0,\varepsilon,s)+\varepsilon) =
\Phi^c(r_0;r_0,L,(\varphi_1^*,h_1^*,\mu^*)(L,\varphi_0,\varepsilon,s),\varepsilon) + \mathcal{O}(\mathrm{e}^{-\eta L}).
\end{equation}
Using (\ref{MatchingSoln}) and Lemma~\ref{lem:CutInvEqn}(i), we see that
\[
\begin{split}
\Phi^c(r_0;r_0,L,(\varphi_1^*,h_1^*,\mu^*)(L,\varphi_0,\varepsilon,s),\varepsilon)
&= \Phi^c(r_0;r_0,L,\tilde{\varphi}(s),0,\tilde{\mu}(s),\varepsilon)
+ \mathcal{O}((\mathrm{e}^{-\eta L}+\varepsilon)\varepsilon\ln L) \\
&= \tilde{\varphi}(s) + r_0 - L + g^c(L,\tilde{\varphi}(s),0,\tilde{\mu}(s),\varepsilon) + \mathcal{O}(\varepsilon+\varepsilon^2\ln L).
\end{split}
\]
Substituting this expression into (\ref{LSoln}) and using (\ref{MatchingSoln}) to estimate $a_0^*$, we see that (\ref{LSoln}) can be written as
\[
\varphi_0 + r_0 + \mathcal{O}(\mathrm{e}^{-\eta L}+\varepsilon) =
\tilde{\varphi}(s) + r_0 - L + g^c(L,\tilde{\varphi}(s),0,\tilde{\mu}(s),\varepsilon) + \mathcal{O}(\varepsilon+\varepsilon^2\ln L)
\]
or as
\begin{equation}\label{LSoln2}
\tilde{\varphi}(s) = \varphi_0 + L - g^c(L,\tilde{\varphi}(s),0,\tilde{\mu}(s),\varepsilon) + \mathcal{O}(\varepsilon+\varepsilon^2\ln L+\mathrm{e}^{-\eta L}).
\end{equation}
Rearranging this equation, we obtain
\[
\tilde{\varphi}(s) = \varphi_0 + L \left( 1-\frac{g^c(L,\tilde{\varphi}(s),0,\tilde{\mu}(s),\varepsilon)}{L} + \mathcal{O}\left(\frac{\varepsilon}{L}+\frac{\varepsilon^2\ln L}{L}+\mathrm{e}^{-\eta L}\right)\right).
\]
Since $g^c(L,\varphi,0,\mu,\varepsilon)=\mathcal{O}(\varepsilon\ln L)$ uniformly in $\varphi$, the right-hand side converges to $\infty$ as $L\to\infty$. Since $\tilde{\varphi}(s)\to\infty$ as $s\to\infty$, we can use the intermediate value theorem to conclude that for each sufficiently large $s$ there is an $L=L^*(\varphi_0,\varepsilon,s)$ that satisfies (\ref{LSoln2}), and we have $L^*(\varphi_0,\varepsilon,s)\to\infty$ as $s\to\infty$. Hence, the functions $(a_0^*,b_0^*,h_0^*,\varphi_1^*,\mu^*,h_1^*)(L,\varphi_0,\varepsilon,s)$ evaluated at $L = L^*(\varphi_0,\varepsilon,s)$ satisfy (\ref{r_0Matching1})--(\ref{r_0Matching4}) and (\ref{LMatching2}). To prove the claim about $\tilde{g}_\mathrm{lift}(L,\mu,\varepsilon)$, we solve (\ref{LSoln2}) for $\tilde{\varphi}(s)$ as a function of $L$ (which can be done as above by reversing the role of $\tilde{\varphi}$ and $L$) and substitute the result into the right-hand side of (\ref{LSoln2}). The claims of the theorem now follow.
\end{proof}

Since we solved equation (\ref{r_0Matching1}) using the intermediate value theorem, it is not clear whether $\Gamma_\mathrm{pulse}^{\varphi_0,\varepsilon}$ is a smooth manifold. We note, however, that the bifurcation curves $\Gamma_\mathrm{pulse}^{\varphi_0,\varepsilon}$ are indeed smooth and unique whenever $\varepsilon\ln L$ is sufficiently small, as we can then differentiate the left-hand side of (\ref{LSoln}) with respect to $L$ and obtain a bound for the derivative.


\section{Vector field on the invariant manifold} \label{sec:IntegralManifold}

In this section, we revisit the vector field on the two-dimensional invariant manifold $\mathcal{P}_\mathrm{e}(\mu,\varepsilon)$. Using the coordinates $(v^c,v^h)$ that parametrize $\mathcal{P}_\mathrm{e}(\mu,\varepsilon)$, we recall that the vector field on $\mathcal{P}_\mathrm{e}(\mu,\varepsilon)$ is given by
\begin{eqnarray}\label{InvariantEqn2}
v^c_x & = & 1 + \frac{\varepsilon}{x}h_1^c(v^c,v^h,\mu,\varepsilon) + \frac{\varepsilon^2}{x^2}h_2^c(x,v^c,v^h,\mu,\varepsilon) \\ \nonumber
v^h_x & = & \frac{\varepsilon}{x}h_1^e(v^c,v^h,\mu,\varepsilon) + \frac{\varepsilon^2}{x^2}h_2^e(x,v^c,v^h,\mu,\varepsilon).
\end{eqnarray}
Our goal is to prove the following theorem which shows that this vector field can be transformed into a simpler form using averaging.

\begin{theorem} \label{thm:AveragingEnergy}
For each fixed $r_0>0$, the transformation
\[
v^c = x + w^c + \varepsilon W^c(x,w^c,w^h,\mu), \qquad
v^h = w^h + \varepsilon W^h(x,w^c,w^h,\mu),
\]
where $W^c$ and $W^h$ are defined in (\ref{WFunction}) below, transforms (\ref{InvariantEqn2}) into the system
\begin{equation} \label{AvgInvMan}
\begin{split}
w^c_x &= \frac{\varepsilon}{x}\left[S^c(w^h,\mu) + \varepsilon F^c(x,w^c,w^h,\mu,\varepsilon)\right], \\
w^h_x &= \frac{\varepsilon}{x}\left[S^h(w^h,\mu) + \varepsilon F^h(x,w^c,w^h,\mu,\varepsilon)\right],
\end{split}
\end{equation}
where
\[
S^h(h,\mu) = \frac{1}{p(\mu,h)}\int_0^{p(\mu,h)} \langle \nabla_u \mathcal{H}(\gamma(x,\mu,h)),g(\gamma(x,\mu,h),\mu,0)\rangle \,\mathrm{d}x.
\]
The functions $F^{c,h}$ are smooth, $p(\mu,w^h)$-periodic with respect to $w^c$, and uniformly bounded in $(x,w^h,\mu,\varepsilon)\in[r_0,\infty)\times K_\mathrm{e}\times J\times[0,\varepsilon_0)$.
\end{theorem}

In the remainder of this section, we will prove this theorem. Let $v^c(x)=x+\tilde{v}^c(x)$ so that $v^c_x=1+\tilde{v}^c_x$, and (\ref{InvariantEqn2}) becomes
\begin{equation}\label{InvariantEqn3}
\begin{split}
\tilde{v}^c_x =&\frac{\varepsilon}{x}h_1^c(x + \tilde{v}^c,v^h,\mu,\varepsilon) + \frac{\varepsilon^2}{x^2}h_2^c(x,x + \tilde{v}^c,v^h,\mu,\varepsilon), \\
v^h_x =&\frac{\varepsilon}{x}h_1^e(x + \tilde{v}^c,v^h,\mu,\varepsilon) + \frac{\varepsilon^2}{x^2}h_2^e(x,x + \tilde{v}^c,v^h,\mu,\varepsilon).
\end{split}
\end{equation}
Lemma~\ref{lem:FenichelCoords} implies that the nonlinearities appearing in (\ref{InvariantEqn3}) are $p(\mu,v^h)$-periodic in $\tilde{v}^c$ and therefore also in $x$.

\begin{lemma} \label{lem:AvgEquiv}
We have
\[
S^h(v^h,\mu) = \frac{1}{p(\mu,v^h)}\int_0^{p(\mu,v^h)} h_1^e(x+\tilde{v}^c,v^h,\mu,0)\,\mathrm{d}x.
\]
\end{lemma}

\begin{proof}
Inspecting the coordinate transformations of Lemma~\ref{lem:FenichelCoords}, we see that
\[
h_1^e(x+\tilde{v}^c,v^h,\mu,0) = \langle \nabla_u \mathcal{H}(\gamma(x+\tilde{v}^c,\mu,v^h)),g(\gamma(x+\tilde{v}^c,\mu,h),\mu,0)\rangle,
\]
and we conclude that
\[
\int_0^{p(\mu,v^h)} h_1^e(x+\tilde{v}^c,v^h,\mu,0)\,\mathrm{d}x = 
\int_0^{p(\mu,v^h)} \langle \nabla_u \mathcal{H}(\gamma(x + \tilde{v}^c,\mu,v^h)), g(\gamma(x + \tilde{v}^c,\mu,v^h),\mu,0)\rangle\,\mathrm{d}x.
\]
Since the integrand on the right-hand side is $p(\mu,v^h)$-periodic in $\tilde{v}^c$, the integral does not depend on $\tilde{v}^c$, and it follows that the right-hand side is equal to $p(\mu,v^h)S^h(v^h,\mu)$.
\end{proof}

Define
\[
S^c(v^h,\mu) = \frac{1}{p(\mu,v^h)}\int_0^{p(\mu,v^h)} h_1^c(x+ \tilde{v}^c,v^h,\mu,0)\,\mathrm{d}x.
\]
The functions
\begin{equation} \label{FTilde}
\tilde{F}^j(x,\tilde{v}^c,v^h,\mu) := h_1^j(x + \tilde{v}^c,v^h,\mu,0) - S^j(v^h,\mu), \qquad j=c,h
\end{equation}
are then $p(\mu,v^h)$-periodic in $x$, and their average over $x$ vanishes. The next result establishes that the coordinate transformations we will employ to derive the averaged vector field and their derivatives are bounded.

\begin{lemma} \label{lem:W}
There exists a $C > 0$ such that for every $r_0 > 0$ the functions
\begin{equation} \label{WFunction}
W^j(x,\tilde{v}^c,v^h,\mu) := \int_{r_0}^x \frac{\tilde{F}^j(s,\tilde{v}^c,v^h,\mu)}{s} \mathrm{d}s, \qquad j=c,h
\end{equation}
and their derivatives with respect to $\tilde{v}^c$ and $v^h$ are bounded by $C/r_0$ for $(x,\tilde{v}^c,v^h,\mu)\in[r_0,\infty)\times\mathbb{R}\times K_\mathrm{e}\times J$.
\end{lemma}

\begin{proof}
Integrating (\ref{WFunction}) by parts gives
\[
W^j(x,\tilde{v}^c,v^h,\mu) = \frac{1}{x}\int_{r_0}^x \tilde{F}^j(s,\tilde{v}^c,v^h,\mu) \,\mathrm{d}s + \int_{r_0}^x \frac{1}{s^2}\int_{r_0}^s \tilde{F}^j(t,\tilde{v}^c,v^h,\mu) \,\mathrm{d}t\,\mathrm{d}s, \qquad j=c,h.
\]
Since the average over $x$ of the functions $\tilde{F}^j(x,\tilde{v}^c,v^h,\mu)$ defined in (\ref{FTilde}) vanishes for each $(\tilde{v}^c,v^h,\mu)$, there exists a constant $C_W>0$ such that the functions
\[
s \longmapsto \int_{r_0}^s \tilde{F}^j(t,\tilde{v}^c,v^h,\mu)\,\mathrm{d}t, \qquad j=c,h
\]
are bounded by $C_W$ for $s\geq r_0$, $\tilde{v}^c\in\mathbb{R}, v^h\in K$ and $\mu\in J$. For each $x\geq r_0$, we therefore have
\[
|W^j(x,\tilde{v}^c,v^h,\mu)| \leq \frac{C_W}{r_0} + C_W \int_{r_0}^x \frac{1}{s^2}\,\mathrm{d}s = \frac{C_W}{r_0}  + \frac{C_W}{r_0} - \frac{C_W}{x} \leq \frac{2C_W}{r_0}.
\]
We can estimate the derivatives of $W^j$ with respect to $(\tilde{v}^c,v^h)$ in the same way as the derivatives are again periodic in $x$ for each fixed $(\tilde{v}^c,v^h)$.
\end{proof}

We can now complete the proof of Theorem~\ref{thm:AveragingEnergy}.

\begin{proof}[Proof of Theorem~\ref{thm:AveragingEnergy}]
We define
\[
W(x,w^c,w^h,\mu) :=
\begin{pmatrix} W^c(x,w^c,w^h,\mu) \\ W^h(x,w^c,w^h,\mu) \end{pmatrix}
\]
and change coordinates according to
\[
\begin{pmatrix} \tilde{v}^c \\ v^h \end{pmatrix} =
\begin{pmatrix} w^c \\ w^h \end{pmatrix} + \varepsilon W(x,w^c,w^h,\mu)
\]
so that
\begin{equation}\label{AvgProof1}
\begin{pmatrix} \tilde{v}^c_x \\ v^h_x \end{pmatrix} =
\left(I + \varepsilon D_{(w^c,w^h)}W(x,w^c,w^h,\mu)\right)
\begin{pmatrix} w^c_x \\ w^h_x \end{pmatrix} + \varepsilon W_x(x,w^c,w^h,\mu),
\end{equation}
where $I$ denotes the $2\times2$ identity matrix. Differentiating (\ref{FTilde}) in $x$ shows that
\[
W^{c,h}_x(x,\tilde{v}^c,v^h,\mu) = \frac{\tilde{F}^{c,h}(x,\tilde{v}^c,v^h,\mu)}{x} = \frac{1}{x}(h_1^{c,e}(x + \tilde{v}^c,v^h,\mu,0) - S^{c,h}(v^h,\mu)).
\]
Using this identity and (\ref{InvariantEqn3}), we can rewrite (\ref{AvgProof1}) to obtain
\begin{equation}\label{AvgProof2}
\begin{split}
(I + \varepsilon & D_{(w^c,w^h)}W(x,w^c,w^h,\mu))
\begin{pmatrix} w^c_x \\ w^h_x \end{pmatrix} =
\begin{pmatrix} \tilde{v}^c_x \\ v^h_x \end{pmatrix} - \varepsilon W_x(x,w^c,w^h,\mu) \\
&= \frac{\varepsilon}{x}\begin{pmatrix} h_1^c(x + w^c + \varepsilon W^c(x,w^c,w^h,\mu),w^h + \varepsilon W^h(x,w^c,w^h,\mu),\mu,\varepsilon) \\
h_1^e(x + w^c + \varepsilon W^c(x,w^c,w^h,\mu),w^h + \varepsilon W^h(x,w^c,w^h,\mu),\mu,\varepsilon) \end{pmatrix} \\
&\quad + \frac{\varepsilon^2}{x^2}\begin{pmatrix} h_2^c(x,x + w^c + \varepsilon W^c(x,w^c,w^h,\mu),w^h + \varepsilon W^h(x,w^c,w^h,\mu),\mu,\varepsilon) \\
h_2^e(x,x + w^c + \varepsilon W^c(x,w^c,w^h,\mu),w^h + \varepsilon W^h(x,w^c,w^h,\mu),\mu,\varepsilon) \end{pmatrix} \\
&\quad - \frac{\varepsilon}{x}
\begin{pmatrix} h_1^c(x + w^c,w^h,\mu,0) \\ h_1^e(x + w^c,w^h,\mu,0) \end{pmatrix}
+ \frac{\varepsilon}{x}\begin{pmatrix} S^c(w^h,\mu) \\ S^h(w^h,\mu) \end{pmatrix} \\
&=: \frac{\varepsilon}{x}
\begin{pmatrix}
S^c(w^h,\mu) + \varepsilon\tilde{F}^c_1(x,w^c,w^h,\mu,\varepsilon)\\
S^h(w^h,\mu) + \varepsilon\tilde{F}^h_1(x,w^c,w^h,\mu,\varepsilon)\end{pmatrix}.
\end{split}
\end{equation}
Lemma~\ref{lem:W} implies that $D_{(w^c,w^h)}W(x,w^c,w^h,\mu)$ is bounded for $x\geq r_0$, and $I+\varepsilon D_{(w^c,w^h)}W(x,w^c,w^h,\mu)$ is therefore invertible for all $(x,w^c,w^h,\mu)$ and all sufficiently small $\varepsilon\geq0$. Thus, (\ref{AvgProof2}) is equivalent to
\[
\begin{split}
\begin{pmatrix} w^c_x \\ w^h_x \end{pmatrix} &=
\frac{\varepsilon}{x} \left(I + \varepsilon D_{(w^c,w^h)}W(x,w^c,w^h,\mu)\right)^{-1}
\begin{pmatrix}
S^c(w^h,\mu) + \varepsilon\tilde{F}^c_1(x,w^c,w^h,\mu,\varepsilon) \\
S^h(w^h,\mu) + \varepsilon\tilde{F}^h_1(x,w^c,w^h,\mu,\varepsilon) \end{pmatrix} \\
&= \frac{\varepsilon}{x}
\begin{pmatrix} S^c(w^h,\mu) \\ S^h(w^h,\mu) \end{pmatrix}
+ \frac{\varepsilon}{x}\left[(I+\varepsilon D_{(w^c,w^h)}W(x,w^c,w^h,\mu))^{-1}-I\right]
\begin{pmatrix} S^c(w^h,\mu) \\ S^h(w^h,\mu) \end{pmatrix} \\
& \quad + \frac{\varepsilon^2}{x} \left(I+\varepsilon D_{(w^c,w^h)}W(x,w^c,w^h,\mu)\right)^{-1}
\begin{pmatrix} \tilde{F}^c_1(x,w^c,w^h,\mu,\varepsilon) \\
\tilde{F}^h_1(x,w^c,w^h,\mu,\varepsilon) \end{pmatrix}.
\end{split}
\]
Observing that $(I+\varepsilon D_{(w^c,w^h)}W(x,w^c,w^h,\mu))^{-1}-I=\mathcal{O}(\varepsilon)$ shows that this system is of the form (\ref{AvgInvMan}).
\end{proof}


\section{Persistent versus collapsed snaking} \label{sec:persistent}

In this section, we investigate when radial pulses with plateau length $L$ exist uniformly in $L\gg1$ and $0\leq\varepsilon\ll1$ and when radial pulses fail to exist for all sufficiently large values $L$ for all but a few distinct values of the parameter $\mu$. In particular, we will complete the proofs of Theorem~\ref{thm:collapsedsnaking} and Theorem~\ref{thm:epsilonlarge}. As shown in Theorem~\ref{thm:epsilonsmall2}, snaking will persist uniformly in $L\gg1$ for all $\mu$ provided that all solutions $v(x)$ of (\ref{InvariantEqn2}) with $|v^h(L)|\ll1$ stay in $K$ for $r_0\leq x\leq L$ regardless of the values of $L\gg1$ and $0\leq\varepsilon\ll1$. Using Theorem~\ref{thm:AveragingEnergy} to write $v=w+\varepsilon W(x,w,\mu)$, it suffices to show that solutions $w(x)$ of
\begin{equation} \label{e4}
\begin{split}
w^c_x &= \frac{\varepsilon}{x}\left[S^c(w^h,\mu) + \mathrm{O}(\varepsilon) \right] \\
w^h_x &= \frac{\varepsilon}{x}\left[S^h(w^h,\mu) + \mathrm{O}(\varepsilon) \right]
\end{split}
\end{equation}
with $w^h(L)=0$ satisfy $w^h(x)\in K$ for $r_0\leq x\leq L$. Hence, we will analyse (\ref{e4}) to identify conditions that guarantee this property as well as conditions that guarantee that such solutions will leave the set $K$ in finite time.

\subsection{Persistent versus non-persistent snaking}

\begin{figure}
\centering
\includegraphics{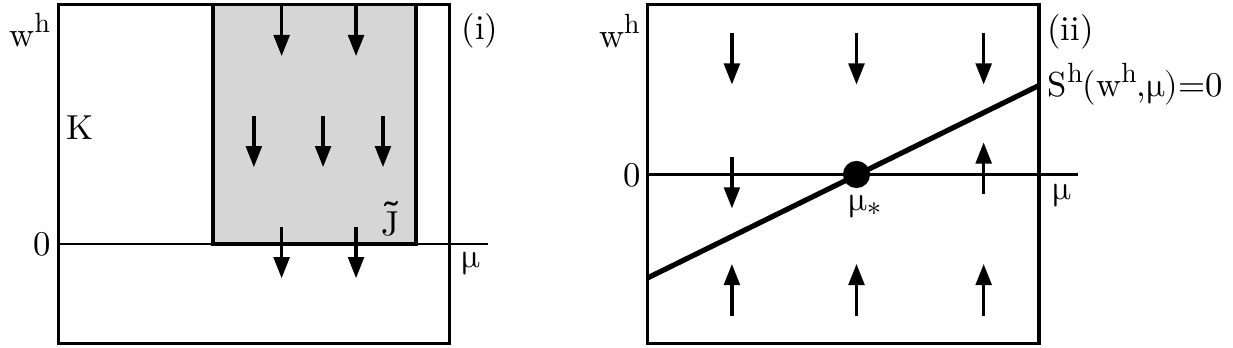}
\caption{Panel~(i) illustrates the hypothesis made in Lemma~\ref{l:nonpersistent} about the existence of an interval $\tilde{J}\subset J$ so that the vector field $S^h(w^h,\mu)$ restricted to $(K\cap\mathbb{R}^+)\times\tilde{J}$ is strictly negative: solutions $w^h(x)$ of (\ref{e4}) with $w^h(L)=0$ cannot remain in $K$ for $r_0\leq x\leq L$ for $L\gg1$. Panel~(ii) illustrates the hypothesis made in Proposition~\ref{p:collapsed} about the function $S^h(w^h,\mu)$ near the parameter value $\mu_*$.}
\label{f:8}
\end{figure}

First, we give conditions so that radial pulses with plateau length $L$ cannot exist for $L>L_\mathrm{min}(\varepsilon)$. We refer to Figure~\ref{f:8}(i) for an illustration of the hypothesis on $S^h$ made in the following lemma.

\begin{lemma} \label{l:nonpersistent}
Assume that there is a closed interval $\tilde{J}\subset J$ so that $S^h|_{(K\cap\mathbb{R}^+)\times\tilde{J}}<0$ or $S^h|_{(K\cap\mathbb{R}^+)\times\tilde{J}}>0$, then for each $0<r_0\ll1$ there are constants $\delta,\varepsilon_0>0$ and a function $L_\mathrm{min}(\varepsilon)$ so that for each function $w(x)$ that satisfies (\ref{e4}) for $(\mu,\varepsilon)\in\tilde{J}\times(0,\varepsilon_0)$ and $L\geq L_\mathrm{min}(\varepsilon)$ with $|w^h(L)|<\delta$ there is a $y\in[r_0,L]$ with $w^h(y)\in K_\mathrm{e}\setminus K$.
\end{lemma}

Note that Theorem~\ref{thm:collapsedsnaking}(i) follows from Theorem~\ref{thm:epsilonsmall2} and Lemma~\ref{l:nonpersistent}.

\begin{proof}
We focus on the case that $S^h|_{(K\cap\mathbb{R}^+)\times\tilde{J}}<0$ as the other case is analogous. In particular, there is a constant $b>0$ so that the right-hand side of the differential equation (\ref{e4}) for $w^h$ satifies
\[
S^h(w^h,\mu) + \mathcal{O}(\varepsilon) \leq -b, \qquad \forall (w^h,\mu,\varepsilon)\in (K\cap\mathbb{R}^+)\times\tilde{J}\times[0,\varepsilon_0].
\]
It follows from Hypothesis~\ref{hyp:Periodic} that we can write $K=[k_-,k_+]$ for some $k_-<0<k_+$. We argue by contradiction and assume that for all $L\gg1$ and $0<\varepsilon<\varepsilon_0$, we have $w^h(x)\leq k_+$ for all $x\in[r_0,L]$. In particular,
\[
k_+ \geq w^h(r_0) = w^h(L) - \int_{r_0}^L \frac{\varepsilon}{y} [S^h(w^h(z),\mu)+\mathcal{O}(\varepsilon)]\,\mathrm{d}z
\geq -\delta + \varepsilon \int_{r_0}^L \frac{b}{y}\,\mathrm{d}z
\geq \varepsilon b \ln \frac{L}{r_0} - \delta.
\]
Since the right-hand side becomes arbitrarily large as $L\to\infty$ for each fixed $\varepsilon>0$, we reach a contradiction. In particular, we conclude that radial pulses with plateau length $L$ cannot exist when
\[
L \geq L_\mathrm{min}(\varepsilon) := r_0 \exp\left( \frac{k_++\delta}{\varepsilon b} \right),
\]
completing the proof of the lemma.
\end{proof}

Next, we will give conditions under which radial pulses with plateau length $L$ exist for all $\mu\in J$, $0\leq\varepsilon\ll1$, and $L\gg1$.

\begin{lemma} \label{l:persistent}
Assume that there is a $w^h_0\in K$ so that one of the following cases is true for all $\mu\in J$:
\begin{compactenum}
\item $S^h(0,\mu)<0$, $S^h(w^h_0,\mu)>0$, and $w^h_0>0$,
\item $S^h(0,\mu)>0$, $S^h(w^h_0,\mu)<0$, and $w^h_0<0$,
\end{compactenum}
then for each $0<r_0\ll1$ there are constants $\delta,\varepsilon_0,L_0>0$ so that solutions $w(x)$ of (\ref{e4}) with $0\leq\varepsilon<\varepsilon_0$ and $|w^h(L)|<\delta$ for $L\geq L_0$ satisfy $w^h(x)\in K$ for all $r_0\leq x\leq L$.
\end{lemma}

\begin{proof}
The claims follow immediately from continuity of $S^h(w^h,\mu)$ and smallness of $\varepsilon$.
\end{proof}

\subsection{Collapsed snaking}

In this section, we will prove Theorem~\ref{thm:collapsedsnaking}(ii), which we restate in the following proposition.

\begin{proposition} \label{p:collapsed}
Assume that Hypotheses~\ref{hyp:Reverser}-\ref{hyp:Pert} are met and that there is a $(\varphi_*,\mu_*)\in\Gamma$ such that $S^h(0,\mu_*)=0$, $S^h_h(0,\mu_*)<0$, $S^h_\mu(0,\mu_*)>0$, and $G_\varphi(\varphi_*,0,\mu_*)\neq0$, then there exists an $\varepsilon_0>0$ such that the following is true for each $0<\varepsilon<\varepsilon_0$: there exists a sequence $(L_m,\mu_m)$ with $L_m\to\infty$ as $m\to\infty$ and $\mu_m$ near $\mu_*$ for all $m$ so that (\ref{ODEPert}) with $\mu=\mu_m$ has a radial pulse with plateau length $L_m$.
\end{proposition}

The hypothesis on $S^h(w^h,\mu)$ made in the proposition, which we illustrate in Figure~\ref{f:8}(ii), implies that there is a smooth function $w^h_*(\mu)$ defined for $\mu\in U_\delta(\mu_*)$ so that $S^h(w^h,\mu)=0$ for $(w^h,\mu)$ near $(0,\mu_*)$ if and only if $w^h=w^h_*(\mu)$. Note that we have $D_\mu w^h_*(\mu_*)>0$. Next, we show that these roots persist as invariant manifolds for the nonautonomous equation (\ref{e4}) for $0<\varepsilon\ll1$.

\begin{lemma} \label{lem:AvgIntegralManifold}
Under the hypotheses of Proposition~\ref{p:collapsed}, there exists a function $W^h_*(x,w^c,\mu,\varepsilon)$ which is smooth in $(w^c,\mu)$, continuous in $(x,\varepsilon)$ for $0\leq\varepsilon\ll1$, and uniformly bounded in all its arguments such that
\[
\{ (w,x) = (w^c,w^h_*(\mu)+\varepsilon W^h_*(\varepsilon \ln x,w^c,\mu,\varepsilon),x):\; (x,w^c)\in[r_0,\infty)\times\mathbb{R} \}
\]
is an invariant manifold of (\ref{e4}) for each $\mu\in U_\delta(\mu_*)$.
\end{lemma}

\begin{proof}
Introducing the coordinates $\tilde{w}^h:=w^h-w^h_*(\mu)$ and $y=\varepsilon\ln x$ so that $x=\mathrm{e}^{y/\varepsilon}$, equation (\ref{e4}) becomes
\[
w^c_x = S^c(w^h_*(\mu)+\tilde{w}^h,\mu) + \mathcal{O}(\varepsilon), \qquad
\tilde{w}^h_x = S^h_\mu(w^h_*(\mu),\mu) \tilde{w}^h + \mathcal{O}(|\tilde{w}^h|^2+\varepsilon).
\]
Since $S^h_\mu(w^h_*(\mu),\mu)<0$ is bounded away from zero uniformly in $\mu\in U_\delta(\mu_*)$, we can apply \cite[Theorem~VII.2.2]{Hale} to conclude that there exists a uniformly bounded function $W^h_*(y,w^c,\mu,\varepsilon)$ that is smooth in $(w^c,\mu)$ and continuous in $(y,\varepsilon)$ such that $\{(w,x)=(w^c,w^h_*(\mu)+\varepsilon W^h_*(\varepsilon \ln x,w^c,\mu,\varepsilon),x)\}$ is an invariant manifold of (\ref{e4}).
\end{proof}

\begin{corollary}\label{cor:InvariantEqn2}
Under the hypotheses of Proposition~\ref{p:collapsed}, there exists $\varepsilon_0>0$ such that for each $(L,\varphi,\mu,\varepsilon)\in[1,\infty)\times\mathbb{R}\times U_\delta(0)\times[0,\varepsilon_0)$ the solution $(v^c_*,v^h_*)(x)\in\mathbb{R}\times K_\mathrm{e}$ of (\ref{InvariantEqn2}) with
\[
v^c_*(L) = \varphi, \qquad
v^h_*(L) = w^h_*(\mu) + \varepsilon W^h_*(L,\varphi,\mu,\varepsilon) + \varepsilon W^h(L,\varphi,w^h_*(\mu)+\varepsilon W^h_*(L,\varphi,\mu,\varepsilon),\mu)
\]
satisfies $v^h_*(x)\in K$ for all $x\in[r_0,L]$.
\end{corollary}

\begin{proof}
Restricting (\ref{AvgInvMan}) to the integral manifold $w^h=w^h_*(\mu)+\varepsilon W^h_*(x,w^c,\mu,\varepsilon)$, equation (\ref{AvgInvMan}) is reduced to
\[
w^c_x = \frac{\varepsilon}{x}\left[S^c(w^h_*(\mu) + \varepsilon W^h_*(x,w^c,\mu,\varepsilon),\mu) + F^c(x,w^c,w^h_*(\mu) + \varepsilon W^h_*(x,w^c,\mu,\varepsilon),\mu,\varepsilon)\right].
\]
For each $\tilde{\varphi}\in\mathbb{R}$ and $L\gg1$, this equation has a unique solution $w^c(x,\tilde{\varphi})$ on $[0,r_0]$ that satisfies $w^c(L,\tilde{\varphi})=\tilde{\varphi}$. Transforming this solution back to the $v^c$-variable using Theorem~\ref{thm:AveragingEnergy}, we find that the corresponding solution $v^c(x,\tilde{\varphi})$ satisfies
\[
v^c(L,\tilde{\varphi}) = L + \tilde{\varphi} + \varepsilon W^c(L,\tilde{\varphi},w^h_*(\mu)+\varepsilon W^h_*(L,\tilde{\varphi},\mu,\varepsilon),\mu).
\]
Using smoothness in $\tilde{\varphi}$, there exists a unique $\tilde{\varphi}$ so that $v^c(L,\tilde{\varphi})=\varphi$. Since $w^h(x)$ stays close to zero for all $x$ by construction, so does $v^h(x)$, and we therefore have $v^h(x)\in K$ for all $x$, which completes the proof.
\end{proof}

\begin{proof}[Proof of Proposition~\ref{p:collapsed}]
We will proceed as in the proof of Theorem~\ref{thm:epsilonsmall2} and therefore focus only on the necessary adjustments. First, it follows from our assumptions that there is a function $(\tilde{\varphi},\tilde{\mu})(s)\in\Gamma$ with $\tilde{\mu}_s(0)\neq0$ so that $(\tilde{\varphi},\tilde{\mu})(0)=(\varphi_*,\mu_*)$. Next, for each $\varphi_1\in\mathbb{R}$, we denote by $(\Phi^c,\Phi^h)(x;r_0,L,\varphi_1,\mu,\varepsilon)$ the solution $(v^c_*,v^h_*)(x)$ constructed in Corollary~\ref{cor:InvariantEqn2} so that
\[
\Phi^c(L;r_0,L,\varphi_1,\mu,\varepsilon) = \varphi_1, \qquad
\Phi^h(L;r_0,L,\varphi_1,\mu,\varepsilon) = w^h_*(\mu) + \mathcal{O}(\varepsilon).
\]
Matching solutions at $x=r_0$ proceeds then as in the proof of Theorem~\ref{thm:epsilonsmall2}, and it remains to solve for the phase and match at $x=L$. Similar to (\ref{LMatching2}), matching at $x=L$ leads to the system 
\begin{eqnarray}
\varphi_1 - \tilde{\varphi}(s) - 2\pi m + \mathcal{O}(\mathrm{e}^{-\eta L}+\varepsilon) & = & 0 \nonumber \\ \label{LMatching3}
\mu - \tilde{\mu}(s) + \mathcal{O}(\mathrm{e}^{-\eta L}+\varepsilon) & = & 0 \\ \nonumber
w^h_*(\mu) + \mathcal{O}(\mathrm{e}^{-\eta L}+\varepsilon) & = & 0,
\end{eqnarray}
where $(\tilde{\varphi}(s)+2\pi m,\tilde{\mu}(s))\in\Gamma_\mathrm{lift}$ for each $m\in\mathbb{N}$. The main difference to the proof of Theorem~\ref{thm:epsilonsmall2} is that the additional free variable $h_1$ that we used to solve (\ref{LMatching2}) is no longer available. Instead we solve (\ref{LMatching3}) for $(\varphi_1,\mu,s)$ near $(\varphi_*+2\pi m,\mu_*,0)$ , which is possible since the Jacobian of the left-hand side of (\ref{LMatching3}) with respect to $(\varphi_1,\mu,s)$, which is given by
\[
\begin{pmatrix} 1 & 0 & -\tilde{\varphi}_s(0) \\ 0 & 1 & -\tilde{\mu}_s(0) \\ 0 & D_\mu w^h_*(\mu_*) & 0 \end{pmatrix} + \mathcal{O}(\mathrm{e}^{-\eta L}+\varepsilon),
\]
is invertible as $D_\mu w^h_*(\mu_*)\neq0$ and $\tilde{\mu}_s(0)\neq0$. Hence, for each $m\in\mathbb{N}$, $L\gg1$, and $0\leq\varepsilon\ll1$, we obtain a solution of (\ref{LMatching3}) in the form
\begin{equation}\label{e10}
(\varphi_1^*,\mu^*,s^*)(L,\varepsilon,m) = (\varphi_*+2\pi m,\mu_*,0) + \mathcal{O}(\mathrm{e}^{-\eta L}+\varepsilon).
\end{equation}
It remains to match the phase at $x=r_0$, which as in the proof of Theorem~\ref{thm:epsilonsmall2} leads to the equation
\[
\varphi_* + 2\pi m - L + \mathcal{O}(\varepsilon\ln L+\mathrm{e}^{-\eta L}) = 0,
\]
which we can solve for $L=L_m$ for each sufficiently large integer $m$ using again the intermediate-value theorem. In particular, $L_m\to\infty$ as $m\to\infty$, and $\mu_m:=\mu^*(L_m,\varepsilon,m)$ is $\mathcal{O}(\mathrm{e}^{-\eta L_m}+\varepsilon)$-close to $\mu_*$.
\end{proof}

\subsection{Extensions to asymptotic radial pulses}

In this section, we prove Theorem~\ref{thm:epsilonlarge}. Thus, we fix $n:=1+\varepsilon>1$ not necessarily close to one and, omitting the dependence of $g$ on $n$, write the differential equation (\ref{ODEPert}) as
\begin{equation}\label{e2}
u_x = f(u,\mu) + \frac{n-1}{x} g(u,\mu).
\end{equation}
We focus on the case that $x$ is large. For each $R\gg1$, we therefore set $x=R+y$ with $y\geq0$ and introduce the new small parameter $\varepsilon$ via $\varepsilon:=1/R\ll1$. Equation (\ref{e2}) then becomes
\begin{equation}\label{e3}
u_y = f(u,\mu) + \frac{\varepsilon(n-1)}{1+\varepsilon y} g(u,\mu), \qquad y\geq0,
\end{equation}
which is of a form similar to (\ref{ODEPert}) except that the denominator $x$ is now replaced by $1+\varepsilon y$. It is straightforward to check that the results for (\ref{ODEPert}) with $x\geq r_0$ stated in \S\ref{sec:Fenichel}, \S\ref{sec:StableManifold}, and \S\ref{sec:IntegralManifold} remain valid for (\ref{e3}) with $y\geq0$. In particular, $1/\varepsilon$-asymptotic radial pulses with plateau length $L$ can exist only if the solutions on the invariant manifold $\mathcal{P}(\mu,\varepsilon)$ belonging to (\ref{e3}) are contained in $K$ for $0\leq y\leq L$.

Under the hypotheses of Theorem~\ref{thm:epsilonlarge}(i), Lemma~\ref{l:nonpersistent} shows that there are constants $L_0,\varepsilon_0>0$ so that for each solution $w(x)$ on the invariant manifold $\mathcal{P}(\mu,\varepsilon_0)$ with $|w^h(L)|\leq\delta$ and $L\geq L_0$ there is a $y_0\in[0,L]$ so that $w^h(y_0)\notin K$. Setting $R_0=1/\varepsilon_0$, this shows that $R_0$-asymptotic radial pulses with plateau length $L$ cannot exist for $L\geq L_0$, completing the proof of Theorem~\ref{thm:epsilonlarge}(i).

Next, we consider Theorem~\ref{thm:epsilonlarge}(ii). Note that we did not impose any conditions on the asymptotic radial pulses at $x=R_*$ (corresponding to $y=0$). Hence, the only requirements are that the underlying solution $(v^c,v^h)(y)$ stays in $K$ for $0\leq y\leq L$ and that the solution $v(y)$ satisfies $v(L)\in W^s_L(0,\mu)$. We can therefore proceed as in the proof of Proposition~\ref{p:collapsed} to solve the matching condition $v(L)\in W^s_L(0,\mu)$ at $y=L$ in the form (\ref{e10}). In contrast to the situation in Proposition~\ref{p:collapsed}, we do not need to match the phase at $y=0$, and asymptotic radial pulses therefore exist for all sufficiently large values of $L$. This completes the proof of Theorem~\ref{thm:epsilonlarge}(ii).


\section{Discussion}

Though our theoretical results apply to a broad class of systems, we focus our discussion on the radial Swift--Hohenberg equation posed on $\mathbb{R}^n$. Amongst our theoretical findings is the proof that the flow on the integral manifold that continues the manifold of roll patterns $U_\mathrm{per}$ to $n>1$ is, to leading order, determined by the vector field $S(h,\mu)$ that is explicitly related to the PDE energy $\mathcal{E}$ and value $h$ of the Hamiltonian $\mathcal{H}$ via
\[
S(h,\mu) = \mathcal{E}(U_\mathrm{per}(\cdot,\mu,h),\mu) - h.
\]
In particular, we showed that $S(0,\mu)$ vanishes precisely at the Maxwell point $\mu=\mu_\mathrm{Max}$. Combining the numerical computation of $S$ together with our theoretical results on the persistence of snaking branches allowed us to conclude that snaking branches for the planar and three-dimensional Swift--Hohenberg equation have to collapse onto the Maxwell point. Our theoretical results also showed that snaking branches persist for all plateau lengths $L\leq\exp(b/|n-1|)$ for $|n-1|\ll1$.

Our analysis does not explain the precise structure of the branches shown in Figure~\ref{f:3} for the planar Swift--Hohenberg equation. In particular, we cannot explain the intermediate stack of isolas nor the fact that the upper snaking branch forms a connected curve. We believe that the specific shape of the snaking diagram is determined by the behavior of solutions away from the invariant manifold of periodic orbits and will therefore likely depend on the global dynamics rather than local properties near the manifold of rolls and the stable manifold of the origin. Investigating the global dynamics away from this manifold would be an interesting project.

We did not investigate the stability of the localized roll solutions for $n>1$. For the Swift--Hohenberg equation in one space dimension, recent work \cite{Makrides} elucidated some of the expected stability properties of localized roll patterns. We expect that the stability results in \cite{Makrides} can be extended to the radial case for $|n-1|\ll1$.

\paragraph{Acknowledgements.}
We are grateful to the referees for constructive comments that helped us strengthen this paper. Bramburger was supported by an NSERC PDF, Altschuler was supported by the NSF through grant DMS-1439786, Avery, Carter, and Sangsawang were supported by the NSF through grant DMS-1148284, and Sandstede was partially supported by the NSF through grants DMS-1408742 and DMS-1714429. 


\end{document}